\documentclass[11 pt,twosided]{amsart}
\usepackage{amsmath,amssymb,mathrsfs,graphicx,bbm}
\usepackage{amsmath}
\usepackage{amsfonts}
\usepackage{latexsym}
\usepackage{amsthm}
\usepackage{amssymb,amscd}
\usepackage{xargs}
\usepackage{tikz}
\usepackage{stmaryrd}
\usepackage{pgf,tikz}
\usepackage{fancybox}
\usepackage{graphicx}
\usepackage{stmaryrd}
\usepackage{color}
\usepackage{enumitem}
\usepackage{subfigure}
\usepackage[colorinlistoftodos]{todonotes}
 \presetkeys{todonotes}%
{inline,backgroundcolor=gray!20,bordercolor=gray!30}{}
\tikzset{/tikz/notestyleraw/.append style={text=black}}

\newtheorem{thm}{Theorem}[section]

\newtheorem{lem}[thm]{Lemma}

\newtheorem{prop}[thm]{Proposition}
\newtheorem{cor}[thm]{Corollary}
\newtheorem{que}[thm]{Open Problem}
\newtheorem{ex}[thm]{Example}
\newtheorem{rmk}[thm]{Remark}
\newcommand{\be}{\begin{eqnarray}}
\newcommand{\ee}{\end{eqnarray}}
\newcommand{\ben}{\begin{eqnarray*}}
\newcommand{\een}{\end{eqnarray*}}
\newcommand{\beal}{\begin{aligned}}
\newcommand{\enal}{\end{aligned}}
\newcommand{\beq}{\begin{equation}}
\newcommand{\eeq}{\end{equation}}

\newcommand{\eps}{\varepsilon}

\newcommand{\lb}{\lambda}
\newcommand{\E}{\mathbb{E}}

\newcommand{\T}{\mathbb{T}}
\newcommand{\R}{\mathbb{R}}

\newcommand{\Z}{\mathbb{Z}}

\newcommand{\om}{\omega}
\newcommand{\dt}{\delta}
\newcommand{\Dt}{\Delta}

\newcommand{\cC}{\mathcal{C}}
\newcommand{\cP}{\mathcal{P}}

\newcommand{\cM}{\mathcal{M}}

\newcommand{\wh}{\widehat }

\usepackage[unicode=true]{hyperref}
\hypersetup{
     colorlinks,
     linkcolor={black!10!blue},
     linkbordercolor = {black!100!blue},
     citecolor={red}
}
\usepackage{mathpazo}

\usepackage{enumitem}
\makeatletter
\def\namedlabel#1#2{\begingroup
    #2%
    \def\@currentlabel{#2}%
    \phantomsection\label{#1}\endgroup
}

\title[Vanishing viscosity limit of HJ equations with nearly optimal discount]{On the vanishing viscosity limit of Hamilton-Jacobi equations with nearly optimal discount}

\thanks{ {\it Statements and Declarations: }The authors declare no competing interests and no data associated.}
\begin{document}
\maketitle


\centerline{Zibo Wang$^*$,\quad Jianlu Zhang$\dagger$}
\medskip
{\footnotesize
\centerline{State Key Laborotary of Mathematical Sciences}
 \centerline{Academy of Mathematics and Systems Science}
 \centerline{Chinese Academy of Sciences, Beijing 100190, China}
  \centerline{{\it Email: }zibowang@amss.ac.cn$^*$,\quad  jellychung1987@gmail.com$\dagger$}  
}

\begin{abstract}
In this paper, we establish the convergence of solutions to the viscous Hamilton-Jacobi equation (with a Tonelli Hamiltonian):
\[
\lb u +H(x, du)=\eps(\lb)\Dt u,\quad \lb>0
\] 
as $\lb\rightarrow 0_+$, once the modulus $\eps(\lb)$ satisfies $\varlimsup_{\lb\rightarrow 0_+}\eps(\lb)/\lb=0$. Such an exponent of $\eps(\lb)$ is nearly optimal in the convergence.\\
\end{abstract}

\section{Introduction}\label{s1}

The study of {\it ergodic constant} was initiated by Lions, Papanicolaou and Varadhan in \cite{LPV}. They proposed a homogenization approach to approximate a constant $c(H)\in\R$, such that the {\it Hamilton-Jacobi equation} 
\beq\label{eq:hj0}
H(x, d_xu)=c(H),\quad x\in \T^n:=\R^n\slash\Z^n.
\eeq
is solvable in the sense of viscosity. Besides, they pointed out that a {\it vanishing discount} scheme can be used to get the ergodic constant. Precisely, for the discounted equation 
\beq\label{eq:dis}
\lb u + H(x, d_x u)=0,\quad x\in \T^n
\eeq
with a parameter $\lb>0$, a unique viscosity solution $u_\lb$ can be obtained due to the {\it comparison principle}. When the Hamiltonian $H:(x, p)\in T^*\T^n\rightarrow\R $ is coercive in $p\in T^*_x\T^n$, they proved that $\lb u_\lb$ uniformly converges to $-c(H)$ as $\lb\rightarrow 0_+$. However, the convergence of  $u_\lb+c(H)/\lb$ as $\lb\rightarrow 0_+$
was just proved recently in \cite{DFIZ}, by using a dynamical approach in light of Aubry-Mather theory. 
After that, there has been a huge activity on this problem and a bunch of important works have been made from different aspects. Among others, we refer to \cite{DW} for a negative discounted version, \cite{DI} for a non-local version, \cite{CDD, IS} for non-compact settings and \cite{CCIZ,C,CFZZ, WYZ, DNYZ} for the nonlinear discounted extensions. Finally, we also mention \cite{IMT,IMT2,MT} for a methodology of general PDE contexts encompassing both viscous and non-viscous cases. 
\medskip

In the current paper, we consider the viscous Hamilton-Jacobi equation 
\beq\label{eq:hj-vis}
\lb u+H(x,d u)=\eps(\lb)\Dt u, \quad\lb>0
\eeq
with a suitable modulus $\eps:\lb\in \R_+\rightarrow\R_+$. Without loss of generality, we can assume $c(H)=0$.  
The asymptotic behaviors of the solution $u_\lb^\eps$ as $\lb\rightarrow 0_+$ will be studied. Such a setting is highly related to the {\it vanishing viscosity} process, which is a long-standing topic on the viscosity solutions. For the Cauchy problem, the convergence of solutions 
in the vanishing viscosity process (as $\eps\rightarrow 0_+$) was given in \cite{H,O,V}, then an $O(\sqrt\eps)-$convergence rate was obtained in \cite{Fle,L,CL,IK,E,T} by using differential games approach, doubling variable technique and nonlinear adjoint method respectively. When the Hamiltonian is convex,  better rates can be expected. When the initial condition is assumed to be semiconcave, a one-sided $O(\eps)-$rate was provided in \cite{L,Cal} and an $O(\eps)-$rate was obtained in the norm $L_t^\infty(L_x^1)$ by \cite{CGM,Gof}. More recently, \cite{CD} and \cite{CG} improved the rate to $O(\eps|\log \eps|)$ by using the large deviation theory and the nonlinear adjoint method respectively. We also refer to \cite{Tu1,Go,HT} for the recent progress on the vanishing viscosity problem in boundary problems.

For the vanishing viscosity process of static Hamilton-Jacobi equations, \cite{Tr} obtained an $O(\eps)-$rate in the sense of average. In a recent work \cite{CG2}, an $O(\eps)-$rate was obtained by transferring the computation in \cite{CG} directly. Both these works need a non-vanishing discount. For equations without discount, the convergence of solutions (by adding additive constants) is generally unknown. This is because a {\it stochastic ergodic constant} will be involved and its asymptotic behavior influences convergence of solutions, see e.g. \cite{G, IM, TZ} and Section \ref{s2} below. Nonetheless, the convergence under specific assumptions was proved in \cite{AIPM, B, JKM}.\medskip

As we can see, the discount term plays a crucial role in the vanishing viscosity process, so we are eager to get the minimal exponent of $\eps(\lb)$ for \eqref{eq:hj-vis}, such that the convergence holds as $\lb\rightarrow 0_+$. Although the ergodic constant in \eqref{eq:hj-vis} is fixed once for all, the stochastic ergodic constant for different values of $\eps$ will be involved during our analysis. Moreover, the limit of solutions as $\lb\rightarrow 0_+$ proves to be a specific solution of \eqref{eq:hj0} closely relevant to the distribution of Mather measures. \medskip

Finally, we would like to refer the readers to \cite{Kr,Lax,Leo,Q,HT} for more applications of the vanishing viscosity limit in different contexts, like the systems of conservative laws, Burgers equations, the Schr\"odinger bridge problem, convergence rate of homogenization and state-constraint optimal controls, etc.


\subsection{Setting and Main result}
Throughout this article, we assume the Hamiltonian $H:(x,p)\in T\T^n\rightarrow\R$ satisfying the followings:
\begin{itemize}
\item [(H1)]\textbf{(Positive definiteness)} For every $(x,p)\in T^*\T^n$, the Hessian matrix $\partial_{pp}H$ is positive definite.
	\item [(H2)] \textbf{(Superlinearity)} For every $x\in\T^n$, $\lim_{|p|\rightarrow+\infty}H(x,p)/|p|=+\infty$, where $|\cdot|$ is the Euclidean norm.
	\item [(H3)] \textbf{(Smoothness)} $H\in C^2(T^*\T^n,\R)$.
	\end{itemize}
Such a kind of Hamiltonian is known as {\it Tonelli Hamiltonians} due to its variational meanings \cite{F}.


\begin{thm}\label{thm:1}
For any Tonelli Hamiltonian $H$, the solution $u_\lb^\eps$ of \eqref{eq:hj-vis} satisfies the following estimate 
\beq\label{eq:qua-est}
\|u^\eps_\lb-u_\lb\|\lesssim \max\big\{\eps/\lb,\eps|\log\eps|\big\}
\eeq
where $u_\lb\in C(\T^n,\R)$ is the unique solution of \eqref{eq:dis} and $\|\cdot\|$ is the sup norm. Consequently, $u_\lb^\eps$ converges to a uniquely identified solution $u$ of \eqref{eq:hj0} as $\lb\rightarrow 0_+$, as long as the modulus $\eps(\lb)$ satisifes $\varlimsup_{\lb\rightarrow 0_+}\eps(\lb)/\lb=0$
 (for brevity we notate it by $\eps(\lb)\leq o(\lb)$). 
\end{thm}
\begin{rmk}
\begin{itemize}
\item The limit $u(x)$ has an expression in \eqref{eq:lim-sol}, which relates with the distribution of Mather measures associated with \eqref{eq:hj0}, see Section \ref{s2} for more details.
\item 
 In Section \ref{s5}, we supplied numerical investigations implying the optimality of the quantitative estimate \eqref{eq:qua-est} (by computing the lower bound of $\|u^\eps_\lb-u_\lb\|$ for a stochastic pendulum). The convergence of $u_\lb$ to $u$ as $\lb\rightarrow 0_+$ was proved in \cite{DFIZ}, whereas the general convergence rate is usually impossible to obtain. We refer the readers to \cite{MS} for the rate of $\|u_\lb-u\|$ in some special examples.  
 \item It's remarkable that by using the method in Theorem 1.47 of \cite{Tr}, we can easily get the estimate of $u_\lb^\eps-u_\lb$ in the sense of average, see Section \ref{s2} for more details.
\item From Section \ref{s2}, we can see that $\varlimsup_{\lb\rightarrow 0_+}\eps(\lb)/\lb<+\infty$, namely $\eps(\lb)\leq O(\lb)$  is necessary to ensure the existence of accumulating functions of $u_\lb^\eps$ as $\lb\rightarrow 0_+$. That's why we say $\eps(\lb)\leq o(\lb)$ is nearly optimal. As a remainder, the following open problem is left to be solved:
\begin{que}
For $\eps(\lb)\leq O(\lb)$, can we prove the convergence of $u_\lb^\eps$ as $\lb\rightarrow 0_+$? If not, can we find a Hamiltonian possessing multiple limits of  $u_\lb^\eps$ as $\lb\rightarrow 0_+$ ? 
\end{que}
\item In our proof of Theorem \ref{thm:1}, the positive definiteness of $\partial_{pp}H$ is essentially required. However, evidence in \cite{CG} implies that for Hamiltonians of the form $H(x,p)=h(p)+f(x)$, the convexity of $h(p)$ can be reduced, which leads to a worse estimate than \eqref{eq:qua-est} then further imposes a more restrictive $\eps(\lb)$ condition.
These facts are beyond the interest of this current work and we present them here to stimulate further discussions. 

\end{itemize}
\end{rmk}

\subsection*{Outline of the paper} 
This paper is organized as the following: In Section \ref{s2}, we review some basic conclusions on the viscosity solutions and Mather measures. In Section \ref{s3}, we give a proof of Theorem \ref{thm:1} by using the sup-convolution approach. In section \ref{s4}, we give a different proof of Theorem \ref{thm:1} by using the adjoint method. We also discuss the difference between these two approaches in this section. In Section \ref{s5},  numerical computations are made for a stochastic pendulum, which shows the optimality of the convergence rate of $\|u_\lb^\eps-u_\lb\|$.

\subsection*{Acknowledgements} 
This work is supported by the National Key R\&D Program of China (No. 2022YFA1007500) and the National Natural Science Foundation of China (No. 12231010, 11901560). ZW is additionally funded by China Postdoctoral Science Foundation (No. 2024M763476) and supported by the National Center of Mathematics and Interdisciplinary Sciences. Both the authors would like to thank Professor Hung Vinh Tran for reminding us of Theorem 1.47 in his textbook \cite{Tr}, which helps improving the exponent in our earlier versions. We also thank Professor Hiroyoshi Mitake for valuable discussions on the rationality of our assumptions, as well as Dr. Jianyu Chen for helpful discussions on the numerical methods.\\

\section{Preliminary: variational properties of viscosity solutions}\label{s2}

Throughout this section, we assume $\eps,\lb\in(0,1]$ are two individual arguments. Recall that 
there exists a unique ergodic constant $c(\eps)\in\R$ such that 
\beq\label{eq:hj-vis-0}
H(x, du)=\eps\Dt u+c(\eps),\quad\eps> 0
\eeq
admits a solution. It's proved in \cite{G}, any two solutions of \eqref{eq:hj-vis-0} differ by a constant. Moreover, any solution of \eqref{eq:hj-vis-0} is $C^2-$smooth. So we can always choose a solution $\om_0^\eps\in C^2(\T^n,\R)$ with $\om_0^\eps(0)=0$. Due to the {\it Perron method} and the {\it comparison principle} (see \cite{CIL,Tr}), there exists a unique solution $\om_\lb^\eps$ of 
\beq\label{eq:hj0-vis}
\lb u+H(x, du)=\eps\Dt u+c(\eps),\quad \eps,\lb>0
\eeq
which satisifes $\om_0^\eps-\|\om_0^\eps\|\leq \om_\lb^\eps\leq \om_0^\eps+\|\om_0^\eps\|$. Moreover, the following two conclusions can be drawn:
\begin{lem}[Theorem 1.1 of \cite{TZ}]\label{lem:c-eps}
$\lim_{\eps\rightarrow 0_+}c(\eps)=c(H)$. Moreover, $c(\eps)$ is uniformly Lipschitz for $\eps\in[0,1]$.
\end{lem}
\begin{lem}[Bernstein's method]\label{lem:bern}
There exists a constant $C_1>0$ uniform w.r.t. $\eps,\lb \in(0,1]$ such that 
\beq\label{eq:1st-d}
\| d\om_\lb^\eps\|\leq C_1<+\infty.
\eeq
\end{lem}
\proof
This result can be proved by adapting Theorem 3.1 of \cite{AT} to $\T^n$. 
\qed\medskip

As an instant corollary of Lemma \ref{lem:c-eps} and Lemma \ref{lem:bern}, $\{\om_\lb^\eps\}$ are uniformly bounded and uniformly Lipschitz w.r.t. $\eps,\lb\in(0,1]$. On the other side, 
\beq
u_\lb^\eps:=\om_\lb^\eps+\frac{c(\eps)-c(H)}{\lb}
\eeq
is exactly the solution of \eqref{eq:hj-vis}, so we have:
\begin{prop}[Boundedness]\label{prop:bound}
For any modulus $\eps(\lb)\leq O(\lb)$, i.e. $\varlimsup_{\lb\rightarrow 0_+}\eps(\lb)/\lb<+\infty$, the family $\{u_\lb^\eps\}$ are uniformly bounded and uniformly Lipschitz w.r.t. $\eps,\lb\in(0,1]$.
\end{prop}
\proof
As $\frac{c(\eps)-c(H)}{\lb}=\frac{c(\eps)-c(H)}{\eps}\cdot\frac\eps\lb$, then due to Lemma \ref{lem:c-eps}, we conclude the uniform boundedness of $\{u_\lb^\eps\}$. Since $du_\lb^\eps=d\om_\lb^\eps$, due to Lemma \ref{lem:bern} we get the uniform Lipschitzness of $\{u_\lb^\eps\}$.
\qed

\subsection{Mather measure} Due to the {\it Legendre transformation}, we can define the {\it Lagrangian} $L:(x,v )\in T\T^n\rightarrow\R$ via:
 \[
 L(x,v):=\max_{p\in T_x^*\T^n}\big(\langle v,p\rangle-H(x,p)\big).
 \]
For any Tonelli Hamiltonian $H$, the associated Lagrangian $L$ can be proved to be $C^2$ smooth, superlinear and positive definite in $v\in T_x\T^n$ for all $x\in\T^n$ (so called  Tonelli Lagrangian). 

Suppose $\cP(T\T^n,\R)$ is the space of all probability measures on $T\T^n$, then we call a measure $\mu\in\cP(T\T^n,\R)$ {\it closed}, if it satisfies:
\begin{itemize} 
\item $\int_{T\T^n}|v|d\mu<+\infty$;
\item $\int_{T\T^n}\langle v,d\phi\rangle d\mu=0$ for all $\phi\in C^1(\T^n,\R)$.
\end{itemize}
Denote by $\cC$ the set of all closed measures, then
\beq\label{eq:min}
\min_{\mu\in \cC}\int_{T\T^n} L(x,v) d\mu=-c(H)=0
\eeq
and the minimizer of \eqref{eq:min} is called a {\it Mather measure}, which is always achievable (see \cite{DFIZ} for instance). The set of all Mather measures is denoted by $\cM$ in our context.\medskip

In \cite{G}, an analogy of the Mather measure was obtained for the equation \eqref{eq:hj-vis} when $\eps>0$. Precisely, we call a measure $\mu\in\cP(T\T^n,\R)$ {\it $\eps-$closed}, if it satisfies:
\begin{itemize} 
\item $\int_{T\T^n}|v|d\mu<+\infty$;
\item $\int_{T\T^n}\big(\langle v,d\phi\rangle-\eps\Dt\phi\big) d\mu=0$ for all $\phi\in C^2(\T^n,\R)$.
\end{itemize}
Denote by $\cC_\eps$ the set of all $\eps-$closed measures, then
\beq\label{eq:min-eps}
\min_{\mu\in \cC_\eps}\int_{T\T^n} L(x,v) d\mu=-c(\eps).
\eeq
Any minimizer of \eqref{eq:min-eps} will be called an {\it $\eps-$Mather measure}, and the set of all such minimizers is denoted by $\cM_\eps$. 
\begin{lem} [Corollary 2.3 \& Lemma 2.6 of \cite{TZ}]\label{lem:eps-mea}
 $\cM_\eps$ is actually a singleton for $\eps>0$. Any accumulating measure of $\mu_\eps\in\cM_\eps$  as $\eps\rightarrow 0_+$ has to be contained in $\cM$.
\end{lem}

\subsection{Adjoint equation}

Previous definition of Mather measures (resp. $\eps-$Mather measures) is equivalently given by using the duality representation in \cite{IMT,IMT2}. Here we propose another approach to get the Mather measures: the {\it nonlinear adjoint method}. This concept was first proposed by Evans in \cite{E}, and then was generalized in \cite{T,Tr,MT}. Precisely, for any $x_0\in\T^n$ fixed, the {\it adjoint equation}
\be\label{eq:ad-app-2}
\lb\theta_\lb^{\eps}-{\rm div}\Big(\partial_p H(x,d_xu_\lb^{\eps})\theta_\lb^{\eps}\Big)=\eps\Dt \theta_\lb^{\eps}+\lb\dt_{x_0}
\ee
admits a continuous solution $\theta_\lb^{\eps}(x)\geq 0$. Moreover, we can prove that $\int_{\T^n}\theta_\lb^\eps(x)dx=1$ and $\theta>0$ for any $x\in\T^n\backslash\{x_0\}$.

\begin{lem}
\begin{itemize}
\item[(i)] $\int_{\T^n}L(x,\partial_pH(x, du_\lb^\eps))\theta_\lb^\eps dx=\lb u_{\lb}^\eps(x_0)$;
\item[(ii)] For any $\psi\in C^2(\T^n,\R)$, there holds 
\[
\int_{\T^n}\Big(\langle d\psi, \partial_p H(x, du_\lb^\eps)\rangle-\eps\Dt \psi\Big)\theta_\lb^\eps dx=\lb\psi(x_0)-\lb\int_{\T^n}\psi\theta_\lb^\eps dx;
\]
\item[(iii)] $\int_{\T^n}|d_x^2 u_\lb^\eps|^2\theta_\lb^\eps dx\leq C_2/\eps$ for some constant $C_2>0$ uniform w.r.t. $\lb,\eps\in(0,1]$.
\end{itemize}
\end{lem}
\proof
Since $u_\lb^\eps$ is $C^2-$smooth, then 
\ben
\langle \partial_pH(x, du_\lb^\eps), d u_\lb^\eps\rangle&=&L(x,\partial_pH(x, du_\lb^\eps))+H(x, du_\lb^\eps)\\
&=&L(x,\partial_pH(x, du_\lb^\eps))+\eps\Dt u_\lb^\eps-\lb u_\lb^\eps.
\een
Integrating both sides w.r.t. $\theta_\lb^\eps(x) dx$, we get item (i). On the other side, 
\ben
& &\int_{\T^n}\Big(\langle d\psi, \partial_p H(x, du_\lb^\eps)\rangle-\eps\Dt \psi\Big)\theta_\lb^\eps dx\\
&=&-\int_{\T^n} {\rm div}\Big(\partial_p H(x,d_xu_\lb^{\eps})\theta_\lb^{\eps}\Big)\psi dx-\eps\int_{\T^n} \Dt\theta_\lb^\eps\psi dx\\
&=&\lb\psi(x_0)-\lb\int_{\T^n}\psi\theta_\lb^\eps dx
\een
so item (ii) follows. Recall that $u_\lb^{\eps}$ is $C^2$ smooth, there holds
\ben
& &\partial_x H(x, d_xu_\lb^{\eps})+\partial_p H(x, d_xu_\lb^{\eps})\cdot d_x^2u_\lb^{\eps}(x)+\lb d_xu_\lb^{\eps}(x)\\
&=&\eps\Dt (d_xu_\lb^{\eps}).
\een
Multiplying previous equality by $d_xu_\lb^{\eps}$ we get 
\be\label{eq:2nd-est}
& &\langle\partial_x H(x, d_xu_\lb^{\eps}), d_xu_\lb^{\eps}\rangle+\langle\partial_p H, d_x\varphi(x)\rangle+2\lb \varphi(x)\\
&=&\eps(\Dt\varphi-|d_x^2u_\lb^{\eps}|^2)\nonumber
\ee
where $\varphi(x):=\dfrac{|d_xu_\lb^{\eps}|^2}{2}$. Due to Lemma \ref{lem:bern}, there exists a constant $C>0$ such that 
\[
|\langle\partial_x H(x, d_xu_\lb^{\eps}), d_xu_\lb^{\eps}\rangle|\leq C.
\]
 Accordingly, \eqref{eq:2nd-est} implies 
\be\label{eq:2nd-est-2}
& &\langle\partial_p H, d_x\varphi(x)\rangle+\lb \varphi(x)
-\eps \Dt\varphi +\eps |d_x^2u_\lb^{\eps}|^2\leq C'.
\ee
Integrating both sides of \eqref{eq:2nd-est-2} by $\theta_\lb^{\eps}(x)dx$ we get 
\ben
\int_{\T^n} \eps |d_x^2u_\lb^{\eps}|^2\theta_\lb^{\eps}(x)dx&\leq& \lb |\varphi(x_0)|+C'
\een
which further indicates 
\beq\label{eq:2nd-ord-est}
\int_{\T^n}  |d_x^2u_\lb^{\eps}|^2\theta_\lb^{\eps}(x)dx\leq \frac{C_2}{\eps}
\eeq
for some constant $C_2>0$ uniform w.r.t. $\lb,\eps \in(0,1]$.
\qed


Next, we will show how to use $\theta_\lb^\eps$ to find Mather measures. We can define a probability measure $\mu_\lb^\eps$ with $\theta_\lb^\eps$ being the density function, i.e.
\be
\int_{T\T^n}\phi(x, v)d\mu_\lb^\eps(x,v)&:=&\int_{\T^n} \phi\Big(x, \partial_p H\big(x,d u_\lb^\eps(x)\big)\Big)\theta_\lb^\eps(x)dx\nonumber\\
&=&\int_{\T^n} \phi\Big(x, \partial_p H\big(x,d \om_\lb^\eps(x)\big)\Big)\theta_\lb^\eps(x)dx,\quad\forall\,\phi\in  C_c(T\T^n,\R).
\ee
There will hold:
\begin{lem}\label{lem:eps-mea-2}
Fix $\eps>0$. Any accumulating measure of $\mu_\lb^{\eps}$ 
as $\lb\rightarrow 0_+$ will be contained in $\cM_\eps$. 
\end{lem}
\proof
Suppose $\mu_\lb^\eps\rightharpoonup\mu$, then we first show that $\mu\in\cC_\eps$. Firstly, due to Lemma \ref{lem:bern},
\[
\int_{T\T^n}|v| d\mu_\lb^\eps=\int_{\T^n}|\partial_pH( x, d\om_\lb^\eps(x))|\theta_\lb^\eps(x)dx\leq C_3<+\infty
\]
for some constant uniform w.r.t. $\eps,\lb\in(0,1]$. That implies $\int_{T\T^n}|v|d\mu\leq C_3<+\infty$. Secondly, for any $\varphi\in C^2(\T^n,\R)$ we have
\ben
\int_{T\T^n}\langle v, \nabla \varphi(x)\rangle -\eps\Dt \varphi(x) d\mu_\lb^\eps(x,v)&=&\int_{\T^n}\Big(\langle v, \nabla \varphi(x)\rangle -\eps\Dt \varphi(x)\Big) \theta_\lb^\eps(x)dx\\
&=&-\int_{\T^n}\Big[{\rm div}\Big(\partial_p H(x,d_xw_\lb^\eps)\theta_\lb^\eps\Big)+\eps\Dt\theta_\lb^\eps\Big]\varphi(x) dx\\
&=&\int_{\T^n}\lb(\dt_{x_0}-\theta_\lb^\eps) \varphi(x)dx\\
&=&\lb\varphi(x_0)-\lb\int_{T\T^n} \varphi(x) d\mu_\lb^\eps(x,v),
\een
then by taking $\lb\rightarrow 0_+$ we get $\int_{T\T^n}\langle v, \nabla \varphi(x)\rangle -\eps\Dt \varphi(x) d\mu^\eps=0$. Consequently, $\mu^\eps\in\cC_\eps$. To prove \eqref{eq:min-eps}, we can see that 
\ben
\int_{T\T^n} L(x,v) d\mu_\lb^\eps(x,v)&=&\int_{T\T^n} \langle \partial_p H(x,\partial_vL), \partial_vL\rangle-H(x,\partial_vL) d\mu_\lb^\eps(x,v)\\
&=&\int_{T\T^n} \Big(\langle \partial_p H(x,d u_\lb^\eps), d u_\lb^\eps\rangle-H(x, du_\lb^\eps) \Big)\theta_\lb^\eps d x\\
&=&\int_{T\T^n} \Big(\langle \partial_p H(x,d \om_\lb^\eps), d \om_\lb^\eps\rangle-H(x, d\om_\lb^\eps) \Big)\theta_\lb^\eps d x\\
&=&\int_{T\T^n} \Big(\langle \partial_p H(x,d \om_\lb^\eps), d \om_\lb^\eps\rangle+\lb w_\lb^\eps-\eps\Dt \om_\lb^\eps \Big)\theta_\lb^\eps d x-c(\eps)\\
&=& \lb \om_\lb^\eps(x_0)-c(\eps)
\een
which tends to  $-c(\eps)$ as $\lb\rightarrow 0_+$ since $\{\om_\lb^\eps\}$ are uniformly bounded for $\lb\in(0,1]$.
\qed
\begin{cor}
For $\eps=\eps(\lb)\leq O(\lb)$, any accumulating measure of $\mu_\lb^{\eps(\lb)}$ 
as $\lb\rightarrow 0_+$ will be contained in $\cM$.
\end{cor}
\proof
Due to Lemma \ref{lem:eps-mea} and Lemma \ref{lem:eps-mea-2} we instantly get the assertion.
\qed

\subsection{Vanishing discount problem} It's proved in \cite{DFIZ}, the function defined by 
\beq\label{eq:var-dis}
u_\lb(x):=\inf_{\substack{\gamma\in C^{ac}((-\infty,0],\T^n)\\\gamma(0)=x}}\int_{-\infty}^{0} e^{\lb t}L(\gamma,\dot\gamma) dt
\eeq
is the unique solution of \eqref{eq:dis}. As $\lb\rightarrow 0_+$, $u_\lb$ converges to a uniquely identified function 
\beq\label{eq:lim-sol}
u(x)=\sup\Big\{\om \text{ is a subsolution of \eqref{eq:hj0}}\Big|\int_{T\T^n}\om(x) d\mu(x,v)\leq 0,\;\;\forall \mu\in\cM \Big\}
\eeq
presenting as a solution of \eqref{eq:hj0}. As a matter of fact, for $\eps,\lb\in(0,1]$ we can quantitatively estimate the discrepancy between $u_\lb^\eps$ and $u_\lb$:
\begin{prop}\label{prop:osc-eps}
There exists a constant $C_4>0$ uniform w.r.t.  $\eps,\lb\in (0,1]$ such that 
\beq\label{eq:osc-eps}
\|u_\lb^\eps-u_\lb\|\leq C_4\sqrt\eps/\lb.
\eeq
\end{prop}
\proof
Taking the derivative of \eqref{eq:hj-vis} w.r.t. $\eps$, we get 
\[
\langle \partial_p H(x,d_xu_\lb^{\eps}),\partial_x\partial_{\eps}u_\lb^{\eps}\rangle +\lb\partial_{\eps} u_\lb^{\eps}=\Dt u_\lb^{\eps}+\eps\Dt(\partial_{\eps} u_\lb^{\eps}).
\]
Consequently, 
\[
\int_{\T^n}\lb\partial_{\eps} u_\lb^{\eps}\theta_\lb^{\eps} dx-\int_{\T^n}\Big(\eps\Dt(\partial_{\eps} u_\lb^{\eps})-\langle \partial_p H,\partial_x\partial_{\eps}u_\lb^{\eps}\rangle\Big)\theta_\lb^{\eps} dx=\int_{\T^n} \Dt u_\lb^{\eps}\theta_\lb^{\eps} dx 
\]
of which the left hand side equals $\lb\partial_{\eps}u_\lb^\eps(x_0)$. Since $x_0\in \T^n$ is freely chosen, we can impose  $|\partial_{\eps} u_\lb^{\eps}(x_0)|=\max_{x\in \T^n}\|\partial_{\eps} u_\lb^{\eps}(x)\|$. Recall that 
\ben
\bigg|\int_{\T^n}\Dt u_\lb^{\eps}\theta_\lb^{\eps} dx\bigg|&\leq&\int_{\T^n}|\Dt u_\lb^{\eps}|\theta_\lb^{\eps} dx\\
&\leq&\sqrt{\int_{\T^n} |d_x^2u_\lb^{\eps}|^2\theta_\lb^{\eps} dx}\cdot\sqrt{\int_{\T^n}\theta_\lb^{\eps} dx}\\
&=&\sqrt{C_2/\eps}
\een
due to \eqref{eq:2nd-ord-est} and the {\it H\"older's inequality}. That further implies  
\[
\|\partial_{\eps} u_\lb^{\eps}(x)\|\leq \frac1\lb \sqrt{C_2/\eps}.
\]
Integrating above inequality w.r.t.  $\eps\in(0,1]$ we get 
\[
\|u_\lb^{\eps}-u_\lb\|=\lim_{\sigma\rightarrow 0_+}\Big\|\int_\sigma^\eps \partial_s u_\lb^s ds\Big\|\leq \lim_{\sigma\rightarrow 0_+}\int_\sigma^\eps \|\partial_s u_\lb^s \|ds\leq \frac{2\sqrt{C_2\eps}}{\lb}
\]
then taking $C_4:=2\sqrt{C_2}$  we get  the assertion.
\qed

As an easy corollary of Proposition \ref{prop:osc-eps}, the following result can be drawn:
\begin{prop}
Suppose $\eps=\eps(\lb)$ is a modulus satisfying $\varlimsup_{\lb\rightarrow0_+}\eps(\lb)/\lb^2=0$ (i.e. $\eps\leq o(\lb^2)$), then $u_\lb^{\eps(\lb)}$ converges to $u$ given by \eqref{eq:lim-sol} as $\lb\rightarrow 0_+$.
\end{prop}
\proof
Due to Proposition \ref{prop:osc-eps}, 
\[
\|u_\lb^{\eps(\lb)}-u\|\leq \|u_\lb^{\eps(\lb)}-u_\lb\|+\|u_\lb-u\|\leq C_4\frac{\sqrt{\eps(\lb)}}\lb+\|u_\lb-u\|
\]
which implies the convergence of $u_\lb^{\eps(\lb)}$ to $u$ as $\lb\rightarrow 0_+$.
\qed
\begin{rmk}
This result is a straightforward corollary of Section 2 of \cite{T}. It gives a rough bound for $\eps(\lb)$ to ensure the convergence of $u_\lb^{\eps}$ as $\lb\rightarrow 0_+$. We emphasize that in Proposition \ref{prop:osc-eps}, the convexity of $H$ in $p$ is unnecessary. However, for Tonelli Hamiltonians, the bound of $\eps(\lb)$ can be further improved later. As a crucial technical point, the adjoint solution $\theta_\lb^\eps$ can be approximated by the solution of the {\bf Fokker-Planck equation}, see \eqref{eq:adj-sta-evo} and Remark \ref{rmk:adj-sta-evo} for related discussions.
\end{rmk}

Before we start to prove the main theorem, we mention that an improvement result can be given in the sense of average. More precisely, for any $x_0\in\T^n$ fixed, consider the following {\it adjoint equation}
\be\label{eq:ad-app-3}
\lb\sigma_\lb^\eps-{\rm div}\Big(\partial_p H(x,d_xu_\lb^\eps)\sigma_\lb^{\eps}\Big)=\eps\Dt \sigma_\lb^\eps+ \lb r(x+x_0).
\ee
We have the following result.

\begin{prop}

(i) There exists a constant $C_5$ uniformly w.r.t. $\lb$ and $\eps$ such that
\[
\int_{\T^n}|d_x^2 u_\lb^\eps(x)|^2\sigma_\lb^\eps(x) dx\le C_5(1+\lb\|d_x r\|_{L^1(\T^n)}).
\]
(ii) For every $x_0\in\T^n$, there exists a constant $C_6:=\sqrt{C_5}$ uniformly w.r.t. $\lb$ and $\eps$ such that
\[
\Big|\int_{\T^n} (u_\lb^\eps(x)-u_\lb(x))r(x+x_0)dx\Big| \le \frac{C_6\eps}{\lb}\Big(1+\lb\|d_x r\|_{L_1(\T^n)}\Big)^{1/2}
\]
\end{prop}

\proof
Differentiating \eqref{eq:hj-vis} twice w.r.t. $x$, we obtain
\ben
\lb z_\lb^\eps+\langle \partial_p H(x,d_xu_\lb^\eps), d_xz_\lb^\eps\rangle+{\rm tr}(d_x^2 u_\lb^\eps\partial_{pp}H(x,d_xu_\lb^\eps)d_x^2 u_\lb^\eps)\\
+2{\rm tr}(\partial_{xp}H(x,d_xu_{\lb}^\eps)d_x^2u_\lb^\eps)+\partial_{xx}H(x,d_xu_\lb^\eps)=\eps\Delta z_\lb^\eps,
\een
where $z_\lb^\eps:=\Delta u_\lb^\eps$.
By the boundedness of $d_xu_\lb^\eps$, there exists a constant $C_7$ such that 
\[
|2{\rm tr}(\partial_{xp}H(x,d_xu_{\lb}^\eps)d_x^2u_\lb^\eps)|+|\partial_{xx}H(x,d_xu_\lb^\eps)|\le \frac{\theta}{2}|d_x^2 u_\lb^\eps|^2+C_7,
\] 
then the uniformly convexity of $H$ implies that
\be\label{eq:ave-est-1}
\lb z_\lb^\eps+\langle \partial_p H(x,d_xu_\lb^\eps), d_xz_\lb^\eps\rangle+\frac{\theta}{2}|d_x^2 u_\lb^\eps|^2 \le \eps\Delta z_\lb^\eps+C_7. 
\ee
Integrating both sides of \eqref{eq:ave-est-1} by $\sigma_\lb^{\eps}(x)dx$ we get 
\ben
\frac{\theta}{2}\int_{\T^n}|d_x^2 u_\lb^\eps(x)|^2\sigma_\lb^{\eps}(x)dx&\le& C_7- \lb \int_{\T^n} \Delta u_\lb^\eps(x) r(x+x_0)dx\\
&=&C_7+\lb \int_{\T^n}  d_x u_\lb^\eps(x) d_xr(x+x_0)dx.
\een
Therefore, the claim (i) obtained from the boundedness of $d_xu_\lb^\eps$.

To prove the claim (ii), we recall that by taking the derivative of \eqref{eq:hj-vis} w.r.t. $\eps$, one gets 
\[
\langle \partial_p H(x,d_xu_\lb^{\eps}),\partial_x\partial_{\eps}u_\lb^{\eps}\rangle +\lb\partial_{\eps} u_\lb^{\eps}=\Dt u_\lb^{\eps}+\eps\Dt(\partial_{\eps} u_\lb^{\eps}).
\]
Consequently, 
\[
\int_{\T^n}\Big(\lb\partial_{\eps} u_\lb^{\eps}-\eps\Dt(\partial_{\eps} u_\lb^{\eps})+\langle \partial_p H,\partial_x\partial_{\eps}u_\lb^{\eps}\rangle\Big)\sigma_\lb^{\eps} dx=\int_{\T^n} \Dt u_\lb^{\eps}\sigma_\lb^{\eps} dx.
\]
Using integration by parts and H\"older inequality, we have
\ben
\lb\Big|\int_{\T^n}\partial_{\eps} u_\lb^{\eps}(x) r(x+x_0)dx\Big|&=&\Big|\int_{\T^n} \Dt u_\lb^{\eps}(x)\sigma_\lb^{\eps}(x) dx\Big| \\
&\le& \sqrt{\int_{\T^n} |\Dt u_\lb^\eps|^2\sigma_\lb^\eps dx}\sqrt{\int_{\T^n} \sigma_\lb^\eps(x)dx}\\
&\le& \sqrt{C_5}(1+\lb\|d_x r\|_{L^1(\T^n)})^{1/2}.
\een
Finally, integrating above inequality w.r.t.  $\eps\in(0,1]$ we get 
\ben
\Big|\int_{\T^n} (u_\lb^\eps(x)-u_\lb(x))r(x+x_0)dx\Big|&=&\lim_{\sigma\to 0_+}\Big|\int_\sigma^\eps\int_{\T^n} \partial_s u_\lb^s(x)r(x+x_0)dxds\Big|\\
&\le&\frac{\sqrt{C_5}\eps}{\lb}(1+\lb\|d_x r\|_{L^1(\T^n)})^{1/2}.
\een
This ends the proof.
\qed

\section{Quantitative estimate on the vanishing viscosity process}\label{s3}

In this section we give a quantitative estimate for $|u_\lb^\eps-u_\lb|$ for $\eps,\lb\in(0,1]$. Similar consideration can be found in \cite{CD,CG} and we make necessary adaptation and generalization. Formally, 
\beq\label{eq:disc-opt}
u_\lb(x)=\inf_{\substack{\gamma\in C^{ac}([0,+\infty),\T^n)\\\gamma(0)=x}}\int_0^{+\infty}e^{-\lb t}\wh L(\gamma,\dot\gamma) dt 
\eeq
and 
\beq\label{eq:disc-sto-opt}
u_\lb^\eps(x)=\inf_{\substack{d\gamma_\om=\mathfrak u dt+\sqrt{2\eps}d B_\om\\\gamma_\om(0)=x}}\E\int_0^{+\infty}e^{-\lb t}\wh L(\gamma_\om,\mathfrak u) dt 
\eeq
where $\wh L(x,v):=L(x,-v)$ is the {\it symmetric Lagrangian} and $\mathfrak u$ could be any admissible progressively measurable control process (see \cite{L} for more details). Our first conclusion is about the uniform semi-concavity of $u_\lb^\eps$, by using the aforementioned optimal control formula:

\begin{lem}
\label{lem:semiconcave}
$\Dt u_\lb^\eps\leq C_\Dt$ for some constant $C_\Dt$ uniform to $\lb>0,\eps\geq 0$.
\end{lem}
\proof
Recall that the infimum of \eqref{eq:disc-sto-opt} could be obtained by the control 
\[
\mathfrak u_*=-\partial_p H(x,du_\lb^\eps),
\]
see \cite{G} for the proof. Suppose $\gamma_\om(t)$ solves the stochastic differential equation:
\[
\left\{
\begin{aligned}
&d\gamma_\om =u_*(\gamma_\om)dt+\sqrt{2\eps} dB_\om(t),\\
&\gamma_\om(0) = x, 
\end{aligned}
\right.
\]
then for any $y\in\R^n$ with $|y|\ll 1$, we conclude that 
\ben
u_\lb^\eps(x+y)&\leq& \E\Big\{\int_0^1e^{-\lb t}\wh L(\gamma_\om+\chi(t), u_*(\gamma_\om)+\dot\chi) dt +\int_1^{+\infty}\wh L(\gamma_\om, u_*(\gamma_\om)) dt\Big\} \\
&=& \E\int_0^1e^{-\lb t}\wh L(\gamma_\om+(1-t)y, u_*(\gamma_\om)-y) dt+ e^{-\lb}u_\lb^\eps(\gamma_\om(1))
\een
where we take $\chi(t):=(1-t)y$ for $t\in[0,1]$ and consequently $\gamma_\om(t)+\chi(t)$ solves the new stochastic differential equation:
\[
\left\{
\begin{aligned}
&d(\gamma_\om+\chi) =(u_*(\gamma_\om)-y)dt+\sqrt{2\eps} dB_\om(t),\\
&(\gamma_\om+\chi)(0) = x+y.
\end{aligned}
\right.
\]
Similarly, we have 
\ben
u_\lb^\eps(x-y)&\leq& \E\Big\{\int_0^1e^{-\lb t}\wh L(\gamma_\om+\varsigma(t), u_*(\gamma_\om)+\dot\varsigma) dt +\int_1^{+\infty}\wh L(\gamma_\om, u_*(\gamma_\om)) dt\Big\} \\
&=& \E\int_0^1e^{-\lb t}\wh L(\gamma_\om-(1-t)y, u_*(\gamma_\om)+y) dt+ e^{-\lb}u_\lb^\eps(\gamma_\om(1))
\een
where  $\varsigma(t):=(t-1)y$ for $t\in[0,1]$. Therefore,
\ben
& &u_\lb^\eps(x+y)+u_\lb^\eps(x-y)-2 u_\lb^\eps(x)\\
&\leq& \E\int_0^1e^{-\lb t}\wh L(\gamma_\om+(1-t)y, u_*(\gamma_\om)-y) dt\\
& &+ \E\int_0^1e^{-\lb t}\wh L(\gamma_\om-(1-t)y, u_*(\gamma_\om)+y) dt\\
& &-2 \E\int_0^1e^{-\lb t}\wh L(\gamma_\om, u_*(\gamma_\om)) dt\\
&\leq &C_\Dt\cdot\sup_{|v|\leq 1+\|\mathfrak u_*\|}\big( \|L_{xx}(x,v)\|+\|L_{xv}(x,v)\|+\|L_{vv}(x,v)\|\big)|y|^2.
\een
for some constant $C_\Dt>0$ uniform w.r.t. $\eps,\lb\in(0,1]$. Here we used the Taylor's expansion of Tonelli Lagrangian $L(x,v)$ around $(\gamma_\om,\mathfrak u_*(\gamma_\om))$ and the uniform Lipschitzness of $\{u_\lb^\eps\}$ in Lemma \ref{lem:bern}.
\qed
\begin{lem}[Upper bound]\label{lem:upperbound}
$u_\lb^\eps\leq u_\lb+\frac\eps\lb C_\Dt$.
\end{lem}
\proof
Recall that $H(x, du_\lb^\eps)+\lb u_\lb^\eps\leq \eps C_\Dt$, so $u_\lb^\eps-\frac\eps\lb C_\Dt$ is a subsolution of \eqref{eq:dis}. By the comparison principle, this subsolution could not be bigger than $u_\lb$, so we conclude the assertion.
\qed
\begin{rmk}
Due to a similar argument as above, we can prove that $u_\lb$ is semi-concave with a uniform linear modulus, i.e. 
\[
u_\lb(x+y)+u_\lb(x-y)-2u_\lb(x)\leq C_8|y|^2, \quad \forall\ x,y\in\T^n
\]
for some constant $C_8>0$ uniform w.r.t. $\lb\in(0,1]$.
\end{rmk}

\subsection{Lower bound of $u_\lb^\eps-u_\lb$} 
Notice that $u_\lb$ is only semiconcave, we can regularize it by the following sup-convolution:
\beq
u_{\lb,\dt}(x):=\sup_{y\in\R^n} \Big(u_\lb(y)-\frac{|y-x|^2}{2\dt}\Big)
\eeq
where $\dt>0$ is an adjustable parameter. Such a regularization has the following properties:

\begin{lem}[Theorem 8.11 of \cite{Cal}]
\label{lem:ulbdt}
There exist constants $C_i>0$ ($i=9,10,11,12,13$) uniform w.r.t. $\lb\in(0,1]$, such that 
\begin{itemize}
\item $u_{\lb,\dt}$ is $C_9-$Lipschitz, $C_{10}-$ semi-concave and $C_{11}/\dt-$ semiconvex. Consequently, $u_{\lb,\dt}$ is $C^{1,1}$ smooth; 
\item $0\leq u_{\lb,\dt}-u_\lb\leq C_{12}\dt$;
\item  $u_{\lb,\dt}$ is a $\dt-$approximated subsolution of \eqref{eq:dis}, i.e.
\beq
\label{eq:hj-dt}
\lb u_{\lb,\dt}+H(x, du_{\lb,\dt})\leq C_{13} \dt.
\eeq
\end{itemize}
\end{lem}
Now we define the function $\vartheta_{\lb,\dt}^\eps:\T^n\rightarrow\R$ by 
\beq
\vartheta_{\lb,\dt}^\eps:=\int_0^1 \partial_p H(x, rd_xu_\lb^\eps+(1-r)d_xu_{\lb,\dt})dr.
\eeq
then consider the stochastic differential equation along the vector field $-\vartheta_{\lb,\dt}^\eps$:
\[
\left\{
\begin{aligned}
&d\gamma_\om(t)=-\vartheta_{\lb,\dt}^\eps(\gamma_\om(t))dt+\sqrt{2\eps}dB_\om(t), \ t\ge 0,\\
&\gamma_\om(0) = x_0.
\end{aligned}
\right.
\]
For fixed $t>0$,  the probability measure induced by this stochastic process $\gamma_\om$ 
is absolutely continuous w.r.t. the Lebesgue measure, i.e.,  $d\mu_{\lb,\dt}^\eps=\rho_{\lb,\dt}^\eps(x,t) dxdt$ where the {\it density function} $\rho_{\lb,\dt}^\eps$ satisfies the following {\it Fokker--Planck equation}:
\beq
\label{eq:fp}
\left\{
\begin{aligned}
&\partial_t \rho_{\lb,\dt}^\eps(x,t)-{\rm div}(\vartheta_{\lb,\dt}^\eps(x)\rho_{\lb,\dt}^\eps(x,t))-\eps\Delta\rho_{\lb,\dt}^\eps(x,t)=0,\\
&\rho_{\lb,\dt}^\eps(x,0)=\delta_{x_0}.
\end{aligned}
\right.
\eeq
The following lemma is based on the upper bounds of the probability density function $\rho_{\lb,\dt}^\eps$ and will be useful later.

\begin{lem}
\label{lem:logrho}
For every $t\ge0$, there exists a constant $C_{14}>0$ uniform w.r.t. $\eps,\lb,\dt\in(0,1]$, such that
\[
\int_{\T^n} |\log{\rho_{\lb,\dt}^\eps(x,t)}|\rho_{\lb,\dt}^\eps(x,t) dx \le C_{14}\left(1+|\log\eps|+|\log t|\right).
\]
\end{lem}
\proof
Recall that by the Lipschitz property of $u_\lb^\eps$ and $u_{\lb,\dt}$, the vector field $\vartheta_{\lb,\dt}^\eps$ is bounded. Hence, by the upper bound estimate of $\rho_{\lb,\dt}^\eps$ (\cite{BKRS}, Corollary 7.2.2), for every $\nu>(n+2)/2$, 
\[
\rho_{\lb,\dt}^\eps \le C\left(1+\frac{1}{\eps}\right)^\nu t^{-n/2}\left(1+\frac{t^{2\nu}}{\eps^\nu}\right).
\]
Therefore,
\[
|\log\rho_{\lb,\dt}^\eps(x,t)|\le C\left(1+|\log \eps|+|\log t|\right).
\]
Then the result follows by that $\rho_{\lb,\dt}^\eps$ is a probability density.
\qed

\begin{lem}
\label{lem:est-dtu}
For every $t\ge1$ and  $\tau\in(0,1)$, there exists a constant $C_{15}>0$ uniform w.r.t. $\eps,\lb,\dt\in(0,1]$ such that the following estimate holds
\[
\int_\tau^t\int_{\T^n} e^{-\lb s}\Delta u_{\lb,\dt}(x)\rho_{\lb,\dt}^\eps(x,s)dxds\ge-C_{15}\left(1+\frac{1}{\lb}+|\log\eps|+|\log\tau|+|\log t|\right).
\]
\end{lem}
\proof
Since the probability density $\rho_{\lb,\dt}^\eps$ satisfies \eqref{eq:fp}, by directly calculating we have
\ben
& & \frac{d}{dt}\int_{\T^n} e^{-\lb t} \log\rho_{\lb,\dt}^\eps(x,t) \rho_{\lb,\dt}^\eps(x,t) dx\\
&=&\int_{\T^n} e^{-\lb t}\Big[-\lb \log\rho_{\lb,\dt}^\eps(x,t) \rho_{\lb,\dt}^\eps(x,t)+\partial_t\rho_{\lb,\dt}^\eps(x,t)(\log\rho_{\lb,\dt}^\eps(x,t)+1)\Big]dx\\
&=&\int_{\T^n} e^{-\lb t}\Big[-\lb \log\rho_{\lb,\dt}^\eps(x,t) \rho_{\lb,\dt}^\eps(x,t)+{\rm div}\big(\vartheta_{\lb,\dt}^\eps(x)\rho_{\lb,\dt}^\eps(x,t)+\nabla\rho_{\lb,\dt}^\eps(x,t)\big)\\
& &(\log\rho_{\lb,\dt}^\eps(x,t)+1)\Big]dx\\
&=&\int_{\T^n} e^{-\lb t}\big[-\lb \log\rho_{\lb,\dt}^\eps(x,t) \rho_{\lb,\dt}^\eps(x,t)-\vartheta_{\lb,\dt}^\eps(x)\cdot\nabla\rho_{\lb,\dt}^\eps(x,t)-\eps\frac{|\nabla\rho_{\lb,\dt}^\eps(x,t)|^2}{\rho_{\lb,\dt}^\eps(x,t)}\big]dx\\
&=&\int_{\T^n} e^{-\lb t}\left[-\lb \log\rho_{\lb,\dt}^\eps(x,t) \rho_{\lb,\dt}^\eps(x,t)+\nabla\cdot\vartheta_{\lb,\dt}^\eps(x) \rho_{\lb,\dt}^\eps(x,t)-\eps\frac{|\nabla\rho_{\lb,\dt}^\eps(x,t)|^2}{\rho_{\lb,\dt}^\eps(x,t)}\right]dx
\een
where we use integrating by parts in the last two equalities. Then we integrate aforementioned equality from $\tau$ to $t$ and use the positivity of $\rho_{\lb,\dt}^\eps$ to obtain that
\be
\label{est}
& &\int_{\T^n}e^{-\lb t}\log\rho_{\lb,\dt}^\eps(x,t)\rho_{\lb,\dt}^\eps(x,t)dx-\int_{\T^n}e^{-\lb \tau}\log\rho_{\lb,\dt}^\eps(x,\tau)\rho_{\lb,\dt}^\eps(x,\tau)dx\\
&\le& \int_\tau^t \int_{\T^n}e^{-\lb s}\left[-\lb \log\rho_{\lb,\dt}^\eps(x,s)\rho_{\lb,\dt}^\eps(x,s)+\nabla\cdot\vartheta_{\lb,\dt}^\eps(x)\rho_{\lb,\dt}^\eps(x,s)\right]dxds\nonumber.
\ee
From the definition of $\vartheta_{\lb,\dt}^\eps$ we get
\ben
\nabla\cdot\vartheta_{\lb,\dt}^\eps(x)&=&\int_0^1 \big(\partial_x\cdot\partial_p\big)H(x,rd_xu_\lb^\eps+(1-r)d_xu_{\lb,\dt})dr\\
&+&\int_0^1{\rm tr}\Big(\partial_{pp}H(x,rd_xu_\lb^\eps+(1-r)d_xu_{\lb,\dt})\cdot\partial_{xx}(ru_\lb^\eps+(1-r)u_{\lb,\dt})\Big)dr
\een
where $\partial_x\cdot\partial_p:=\sum_{i=1}^n\partial_{x_ip_i}$ for brevity.
Due to Lemma \ref{lem:bern} and Lemma \ref{lem:ulbdt}, there exists a constant $C_{16}>0$ such that 
\[
\big(\partial_x\cdot\partial_p\big)H(x,rd_xu_\lb^\eps+(1-r)d_xu_{\lb,\dt})\le C_{16}.
\]
Moreover, due to the positive definiteness of $\partial_{pp}H$, there should exists two constants $\Theta\geq \theta>0$ such that 
\beq\label{eq:convexity}
\theta\mathbb{I}\preceq \partial_{pp}H(x,rd_xu_\lb^\eps+(1-r)d_xu_{\lb,\dt})\preceq \Theta\mathbb{I}.
\eeq
holds for all $x\in\T^n$ and $r\in[0,1]$, where $\mathbb I$ is the identity matrix of the rank $n$.
Therefore, by using the Lemma \ref{lem:matrix} for trace of matrix below, we have
\ben
\nabla\cdot\vartheta_{\lb,\dt}^\eps(x)&\le&C_{16}+\int_0^1{\rm tr}\Big(\partial_{pp}H(x,rd_xu_\lb^\eps+(1-r)d_xu_{\lb,\dt})\cdot\partial_{xx}(ru_\lb^\eps+(1-r)u_{\lb,\dt})\Big)dr\\
&\le& C_{16}+\int_0^1 \theta {\rm tr}(\partial_{xx}(ru_\lb^\eps+(1-r)u_{\lb,\dt}))dr+\max\{C_{10},C_\Dt\}\cdot(\Theta-\theta)n\\
&=& C_{16}+\frac{\theta}{2}\Delta u_\lb^\eps+\frac{\theta}{2}\Delta u_{\lb,\dt}+C_{10}(\Theta-\theta)n
\een
where we assume $C_{10}\geq C_\Dt$ without loss of generality. 
 Hence, by transferring \eqref{est} we get 
\be
\label{est1}
& &\frac{\theta}{2}\int_\tau^t\int_{\T^n} e^{-\lb s}\Delta u_{\lb,\dt}(x)\rho_{\lb,\dt}^\eps(x,s) dxds\\
&\ge& \int_{\T^n}e^{-\lb t}\log\rho_{\lb,\dt}^\eps(x,t)\rho_{\lb,\dt}^\eps(x,t)dx-\int_{\T^n}e^{-\lb \tau}\log\rho_{\lb,\dt}^\eps(x,\tau)\rho_{\lb,\dt}^\eps(x,\tau)dx\nonumber\\
& &+\int_\tau^t\int_{\T^n} \lb e^{-\lb s} \log\rho_{\lb,\dt}^\eps(x,s)\rho_{\lb,\dt}^\eps(x,s)dxds\nonumber\\
& &-\int_\tau^t\int_{\T^n}e^{-\lb s} \Big( \frac{\theta}{2}\Delta u_\lb^\eps(x)+C_7(\Theta-\theta)n+C_{12}\Big)\rho_{\lb,\dt}^\eps(x,s) dxds \nonumber
\ee
For the first two integrals on the right hand side of \eqref{est1}, by Lemma \ref{lem:logrho}, we have
\[
\int_{\T^n}e^{-\lb t}\log\rho_{\lb,\dt}^\eps(x,t)\rho_{\lb,\dt}^\eps(x,t)dx\ge -C_{14}e^{-\lb t}(1+|\log\eps|+|\log t|),
\]
and
\[
-\int_{\T^n}e^{-\lb \tau}\log\rho_{\lb,\dt}^\eps(x,\tau)\rho_{\lb,\dt}^\eps(x,\tau)dx\ge -C_{14}e^{-\lb \tau}(1+|\log\eps|+|\log\tau|).
\]
Moreover, we have the following estimate for the third integral of \eqref{est1}, 
\ben
&&\int_\tau^t\int_{\T^n} \lb e^{-\lb s} \log\rho_{\lb,\dt}^\eps(x,s)\rho_{\lb,\dt}^\eps(x,s)d xds\\
&\ge& -C_{14}\int_\tau^t \lb e^{-\lb s}(1+|\log\eps|+|\log s|)ds\\
&=&C_{14}(1+|\log\eps|)(e^{-\lb t}-e^{-\lb \tau})-C_{14}\int_\tau^t \lb e^{-\lb s}|\log s| ds \\
&\ge& -C_{14}\left((1+|\log\eps|+|\log\tau|)e^{-\lb \tau}+\frac{2}{\lb}e^{-\lb}\right)
\een
where we used
\ben
-\int_\tau^t \lb e^{-\lb s}|\log s|ds&=& \int_\tau^t |\log s| de^{-\lb s}=e^{-\lb s}|\log s|\big|_\tau^t-\int_\tau^t e^{-\lb s}d |\log s|\\
&=&e^{-\lb t}|\log t|-e^{-\lb\tau}|\log\tau|+\int_\tau^1 \frac{e^{-\lb s}}{s}ds-\int_1^t \frac{e^{-\lb s}}{s} ds\\
&\ge& -e^{-\lb \tau}|\log\tau|+\int_\tau^1 e^{-\lb s} ds-\int_1^t e^{-\lb s}ds \\
&\ge&-e^{-\lb \tau}|\log\tau|-\frac{2}{\lb}e^{-\lb}.
\een
Finally, for the last term on the right hand side of \eqref{est}, by  Lemma \ref{lem:semiconcave} we have
\ben
& &-\int_\tau^t\int_{\T^n}e^{-\lb s} \Big( \frac{\theta}{2}\Delta u_\lb^\eps(x)+C_7(\Theta-\theta)n+C_{12}\Big)\rho_{\lb,\dt}^\eps(x,s) dxds\\
&&\ge - \int_\tau^t e^{-\lb s} \Big(\frac{\theta}{2}C_\Delta+C_{10}(\Theta-\theta) n+C_{16}\Big)ds\\
&&=\Big(\frac{\theta}{2}C_\Delta+C_{10}(\Theta-\theta)n+C_{16}\Big)\frac{e^{-\lb t}-e^{-\lb\tau}}{\lb}\\
&&\ge-\Big(\frac{\theta}{2}C_\Delta+C_{10}\Theta n+C_{16}\Big)\frac{1}{\lb}e^{-\lb\tau}.
\een
Combining the previous estimates together, we get
\ben
&&\frac{\theta}{2}\int_\tau^t\int_{\T^n} e^{-\lb s}\Delta u_{\lb,\dt}(x)\rho_{\lb,\dt}^\eps(x,s)ds\\
&\ge& -C_{14}\Big(e^{-\lb t}(1+|\log\eps|+|\log t|)+2e^{-\lb\tau}(1+|\log\eps|+|\log\tau|)\Big)\\
&&-\frac{1}{\lb}\Big(2C_{14} e^{-\lb}+\Big(\frac{\theta}{2}C_\Delta+C_{10}\Theta n+C_{16}\Big)e^{-\lb\tau}\Big),
\een
which implies the result.
\qed
\begin{lem}[matrix inequality]
\label{lem:matrix}
Let $A, B$ be two $n\times n$-dimensional symmetric matrices. If there exists constants $\theta, \Theta, M>0$ such that $\theta\mathbb{I}\preceq  A\preceq \Theta\mathbb{I}$ and $B\le M\mathbb{I}$, then 
\[
{\rm tr}(AB)\le \theta{\rm tr}(B)+(\Theta-\theta)Mn.
\]
\end{lem}
\proof
For the symmetric matrix $A$, there exists an orthogonal matrix $P$ and a diagonal matrix $D$ such that 
\[A=P^TDP,\]
where the diagonal matrix $D:={\rm diag(d_1,...,d_n)}$ with eigenvalues $d_i$ of $A$. Moreover, $\theta\le d_i\le \Theta, i=1,..,n$ since  $\theta\mathbb{I}\preceq  A\preceq \Theta\mathbb{I}$. By the property of trace, we have
\ben
{\rm tr}(AB)&=&{\rm tr}\big(A(B-M\mathbb{I})+M\mathbb{I}\big)={\rm tr}\big(P^T DP(B-M\mathbb{I})\big)+M{\rm tr}(P^TDP)\\
&=&{\rm tr}\big(DP(B-M\mathbb{I})P^T\big)+M{\rm tr}(D)
\een
Note that $B-M\mathbb{I}\le 0$, we obtain
\[
{\rm tr}(AB)\le \theta{\rm tr}\big(P(B-M\mathbb{I})P^T\big)+\Theta Mn=\theta{\rm tr}(B)+(\Theta-\theta)Mn.
\]
and finish the proof.
\qed

\begin{thm}[Lower bound]\label{thm:low-b}
For any Tonelli Hamiltonian $H$, there exists a positive constant $C_{17}$ uniform w.r.t. $\eps,\lb\in(0,1]$ such that
\[
u_\lb^\eps(x)-u_\lb(x) \ge C_{17}\left(-\frac{\eps}{\lb}+\eps\log\eps\right).
\]
In particular, the solution $u_\lb^\eps$ converges to $u_\lb$ (as $\lb\rightarrow 0_+$) if $\eps\leq o(\lb)$. Moreover, the optimal convergence rate will be achieved by $\eps\lesssim e^{-1/\lb}$.
\end{thm}
\proof
As we know, 
\[
u_\lb^\eps-u_\lb=u_\lb^\eps-u_{\lb,\dt}+u_{\lb,\dt}-u_\lb\ge u_\lb^\eps-u_{\lb,\dt}
\]
where $u_{\lb,\dt}-u_\lb$ is positive by Lemma \ref{lem:ulbdt}. It remains to estimate the lower bound of $u_\lb^\eps-u_{\lb,\dt}$. Subtracting \eqref{eq:hj-vis} to \eqref{eq:hj-dt} one obtain that
\be
\label{eq:hj-vis-dt}
& &-\lb (u_\lb^\eps-u_{\lb,\dt})(x)-\vartheta_{\lb,\dt}^\eps(x)\nabla(u_\lb^\eps-u_{\lb,\dt})(x)+\eps\Delta(u_\lb^\eps-u_{\lb,\dt})(x)\\
&\le& -\eps\Delta u_{\lb,\dt}(x) +C_{13}\dt.\nonumber
\ee
Applying \eqref{eq:fp} and \eqref{eq:hj-vis-dt}, we have
\ben
&&\frac{d}{dt}\int_{\T^n}e^{-\lb t}(u_\lb^\eps-u_{\lb,\dt})(x)\rho_{\lb,\dt}^\eps(x,t)dx\\
&=&\int_{\T^n} -\lb e^{-\lb t}(u_\lb^\eps-u_{\lb,\dt})(x)\rho_{\lb,\dt}^\eps(x,t)+e^{-\lb t}(u_\lb^\eps-u_{\lb,\dt})(x)\partial_t\rho_{\lb,\dt}^\eps(x,t)dx\\
&\le&\int_{\T^n} e^{-\lb t}\left[\left(\vartheta_{\lb,\dt}^\eps(x)\cdot\nabla-\eps\Delta\right)(u_\lb^\eps-u_{\lb,\dt})(x)-\eps\Delta u_{\lb,\dt}(x)+C_{10}\dt\right]\rho_{\lb,\dt}^\eps(x,t)\\
&&+e^{-\lb t}\left[(u_\lb^\eps-u_{\lb,\dt})(x)\left({\rm div}(\nabla\vartheta_{\lb,\dt}^\eps(x)\rho_{\lb,\dt}^\eps(x,t))+\eps\Delta\rho_{\lb,\dt}^\eps(x,t)\right)\right]dx\\
&=&\int_{\T^n} e^{-\lb t}(-\eps \Delta u_{\lb,\dt}(x)+C_{13}\dt)\rho_{\lb,\dt}^\eps(x,t) dx.
\een
By integrating it form $0$ to $t$ we obtain that
\be\label{eq:diff-sol}
& &\int_{\T^n}e^{-\lb t}(u_\lb^\eps-u_{\lb,\dt})(x)\rho_{\lb,\dt}^\eps(x,t)dx-(u_\lb^\eps-u_{\lb,\dt})(x)\\
&\le&-\eps\int_0^t\int_{\T^n} e^{-\lb s}\Delta u_{\lb,\dt}(x)\rho_{\lb,\dt}^\eps(x,s)dx ds+C_{13}\dt\int_0^t e^{-\lb s}ds.\nonumber
\ee
Now we split the integral on the right hand side of \eqref{eq:diff-sol} into two parts. For the integral on time interval $(0,\tau)$, applying the semiconvex property of $u_{\lb,\dt}$ to obtain that
\ben
\eps\int_0^\tau\int_{\T^n} e^{-\lb s}\Delta u_{\lb,\dt}(x)\rho_{\lb,\dt}^\eps(x,s) dxds\ge -\eps\frac{C_{11}}{\delta}\int_0^\tau e^{-\lb s}ds\ge -C_{11}\frac{\tau}{\dt}\eps.
\een
For the integral on time interval $(\tau,t)$,  Lemma \ref{lem:est-dtu} implies that
\ben
\eps\int_\tau^t\int_{\T^n} e^{-\lb s}\Delta u_{\lb,\dt}(x)\rho_{\lb,\dt}^\eps(x,s)dxds\ge-C_{15}\left(1+\frac{1}{\lb}+|\log\eps|+|\log\tau|+|\log t|\right)\eps.
\een
Therefore, by the boundness of $u_\lb^\eps$ and $u_{\lb,\dt}$ there exists a constant $C_{18}>0$ such that
\ben
(u_\lb^\eps-u_{\lb,\dt})(x)&\ge& \int_{\T^n}e^{-\lb t}(u_\lb^\eps-u_{\lb,\dt})(x)\rho_{\lb,\dt}^\eps(x,t)dx-C_{13}\dt\int_0^t e^{-\lb s}ds\\
&&-C_{11}\frac{\tau}{\dt}\eps-C_{15}\left(1+\frac{1}{\lb}+|\log\eps|+|\log\tau|+|\log t|\right)\eps\\
&\ge& -C_{18}\Big( \frac{\dt+\eps}{\lb}+\big(1+\frac{\tau}{\dt}+|\log\eps|+|\log\tau|+|\log t|\big)\eps+e^{-\lb t}\Big).
\een
Let $\delta=\tau=\eps$ and $t=\eps^{-1}$, we have
\[
(u_\lb^\eps-u_{\lb,\dt})(x)\ge-3C_{18}\left(\frac{\eps}{\lb}-\eps\log\eps+e^{-\lb/\eps}\right)\ge 6C_{18}\left(-\frac{\eps}{\lb}+\eps\log\eps\right)
\]
and then by taking $C_{17}=6C_{18}$
\[
u_\lb^\eps-u_\lb=u_\lb^\eps-u_{\lb,\dt}+u_{\lb,\dt}-u_\lb\ge u_\lb^\eps-u_{\lb,\dt}\ge  C_{17}\left(-\frac{\eps}{\lb}+\eps\log\eps\right)
\]
which ends the proof.\qed
\bigskip

\noindent{\it Proof of Theorem \ref{thm:1}:}
Recall that by Lemma \ref{lem:upperbound} and Theorem \ref{thm:low-b}, we have
\[
C_{17}\Big(-\frac{\eps(\lb)}{\lb}+\eps\log\eps\Big)\le u_\lb^\eps-u_\lb\le C_\Delta \frac{\eps(\lb)}{\lb}.
\]
Hence, the difference between $u_\lb^{\eps(\lb)}$ and $u$ could be estimated by
\[
\|u_\lb^{\eps(\lb)}-u\|\leq \|u_\lb^{\eps(\lb)}-u_\lb\|+\|u_\lb-u\|\leq \max\{C_{17},C_\Dt\}\Big(\frac{\eps(\lb)}\lb-\eps\log\eps\Big)+\|u_\lb-u\|
\]
which implies that the convergence of $u_\lb^{\eps(\lb)}$ to $u$ as $\lb\rightarrow 0_+$
once the modulus $\eps(\lb)\leq o(\lb)$.
\qed
\medskip

\section{An alternative proof of Theorem \ref{thm:1}: the adjoint method}\label{s4}
As in Proposition \ref{prop:osc-eps}, taking the derivative of \eqref{eq:hj-vis} w.r.t. $\eps$ and letting $v_\lb^\eps:=\partial_\eps u_\lb^\eps$, we get 
\beq
\label{eq:hj-partialu}
\langle \partial_p H(x,d_xu_\lb^{\eps}),\partial_x v_\lb^{\eps}\rangle +\lb v_\lb^{\eps}=\Dt u_\lb^{\eps}+\eps\Dt v_\lb^{\eps}.
\eeq

Consider the following Fokker--Planck equation
\beq
\label{eq:fp-2}
\left\{
\begin{aligned}
&\partial_t \rho_\lb^\eps(x,t)-{\rm div}(\partial_p H(x,d_xu_\lb^\eps(x))\rho_\lb^\eps(x,t))-\eps\Delta\rho_\lb^\eps(x,t)=0,\\
&\rho_\lb^\eps(x,0)=\delta_{x_0},
\end{aligned}
\right.
\eeq
whose weak solution is the probability density corresponding to the following SDE:
\ben
d\gamma_\om(t)=-\partial_p H(x, d_xu_\lb^\eps(\gamma_\om(t)))dt+\sqrt{2\eps}dB_\om(t), \quad \gamma_\om(0)=x_0.
\een

\begin{rmk}\label{rmk:adj-sta-evo}
Actually, the Fokker--Planck equation \eqref{eq:fp-2} is corresponding to the adjoint equation \eqref{eq:ad-app-2} with
\beq\label{eq:adj-sta-evo}
\theta_\lb^\eps(x)=\lim_{T\to\infty}\int_0^T \lb e^{-\lb t}\rho_\lb^\eps(x,t)dt,
\eeq
which can be obtained by a direct calculation. Precisely, multiplying the Fokker--Planck equation \eqref{eq:fp-2} by $e^{-\lb t}$ to get
\[
\partial_t (e^{-\lb t}\rho_\lb^\eps(x,t))=-\lb e^{-\lb t}\rho_\lb^\eps(x,t)+e^{-\lb t}\Big({\rm div}(\partial_p H(x,d_x u_\lb^\eps(x))\rho_\lb^\eps(x,t))+\Dt \rho_\lb^\eps(x,t)\Big),
\]
and integrating it from $0$ to $T$, we obtain that 
\ben
&&e^{-\lb T}\rho_\lb^\eps(x,T)-\rho_\lb^\eps(x,0)\\
&&=\int_0^T \Big(-\lb e^{-\lb t}\rho_\lb^\eps(x,t)+{\rm div}\big(\partial_p H(x,d_x u_\lb^\eps(x)) e^{-\lb t}\rho_\lb^\eps(x,t)\big)+\eps\Delta(e^{-\lb t}\rho_\lb^\eps(x,t))\Big)dt.
\een
Then by letting $T\to+\infty$, the previous equation becomes
\[
\lb(\theta_\lb^\eps(x)-\delta_{x_0})+{\rm div}(\partial_p H(x,d_xu_\lb^\eps(x))\theta_\lb^\eps(x))=\eps\Delta\theta_\lb^\eps(x).
\]
That is to say that $\theta_\lb^\eps$ satisfies the adjoint equation \eqref{eq:ad-app-2}.
\end{rmk}

In the following, we aim to estimate $\| u_\lb^\eps-u_\lb\|$ by using the adjoint method. First, we recall the following property for $\rho_\lb^\eps$ as in  \cite{CG}. 
\begin{lem}
\label{lem:est-Drho}
For every $t\in(0,1)$, there exists a suitably small $\eps_0:=\eps_0(t)>0$ such that for any $\lb\in(0,1]$ and $\eps\in(0,\eps_0]$, the following 
\beq
\label{eq:est-Drho1}
\int_{\T^n}|\partial_x\rho_\lb^\eps(x,t)|dx \le \frac{\sqrt{2C_{20}}}{\eps}(1+2|\log\eps|+|\log t|)^{1/2}.
\eeq
holds for a uniform constant $C_{20}>0$. Furthermore, there exists a constant $C_{21}>0$ such that
\beq
\label{eq:est-Drho2}
\int_{\T^n}|\partial_x\rho_\lb^\eps(x,t)|dx \le \frac{C_{21}}{\eps}(1+|\log\eps|+|\log t|).
\eeq
\end{lem}
\proof
The proof relies on the estimate of the heat kernel. For completeness and convenience, we state the proof as follows. Recall that as in Lemma \ref{lem:logrho}, there exists a constant $C_{14}>0$ uniform w.r.t. $\eps,\lb\in(0,1]$, such that
\beq
\label{eq:logrho}
\int_{\T^n}|\log \rho_\lb^\eps(x,t)|\rho_\lb^\eps(x,t)dx\le C_{14}(1+|\log\eps|+|\log t|).
\eeq
Testing the equation \eqref{eq:fp-2} by $\log\rho_\lb^\eps$ and integrating it on $\T^n$, we obtain 
\ben
&&\frac{d}{dt}\int_{\T^n}\rho_\lb^\eps(x,t)\log\rho_\lb^\eps(x,t)dx+\eps\int_{\T^n}\frac{|\partial_x \rho_\lb^\eps(x,t)|^2}{\rho_\lb^\eps(x,t)}dx\\
&&=-\int_{\T^n} \partial_p H(x,d_xu_\lb^\eps(x))\partial_x\rho_\lb^\eps(x,t) dx.
\een
Hence, by Young's inequality and the Lipschitz property of $u_\lb^\eps$, there exists a constant $C_{19}>0$ such that
\ben
&&\frac{d}{dt}\int_{\T^n}\rho_\lb^\eps(x,t)\log\rho_\lb^\eps(x,t)dx+\frac{\eps}{2}\int_{\T^n}\frac{|\partial_x \rho_\lb^\eps(x,t)|^2}{\rho_\lb^\eps(x,t)}dx\\
&&\le \frac{1}{2\eps}\int_{\T^n} |\partial_p H(x,d_x u_\lb^\eps(x))|^2\rho_\lb^\eps(x,t)dx\le\frac{C_{19}}{\eps}.
\een
For any fixed $t>0$, we can choose a positive $\eps_0\leq t/2$, and an arbitrary interval $[t_0,t_1]$ containing $t$ with $t_1-t_0=\eps\leq\eps_0 $.  
Integrating pervious inequality on intervals $(t_0,t_1)$ and combining it with \eqref{eq:logrho}, 
we get 
\[
\eps\int_{t_0}^{t_1}\int_{\T^n}\frac{|\partial_x \rho_\lb^\eps(x,t)|^2}{\rho_\lb^\eps(x,t)}dxdt\le C_{21}\left(1+|\log\eps|+|\log t_0|+|\log t_1|+\frac{t_1-t_0}{\eps}\right),
\]
where $C_{20}=\max\{2C_{14}, 2C_{19}\}$.
By the mean value theorem, there exists $\hat{t}\in[t_0,t_1]$ such that
\ben
\int_{\T^n}\frac{|\partial_x\rho_\lb^\eps(x,\hat{t})|^2}{\rho_\lb^\eps(x,\hat{t})}dx&=&\frac{1}{t_1-t_0}\int_{t_0}^{t_1}\int_{\T^n} \frac{|\partial_x\rho_\lb^\eps(x,t)|^2}{\rho_\lb^\eps(x,t)}dxdt\\
&\le& \frac{C_{20}}{(t_1-t_0)\eps}\Big(1+|\log\eps|+|\log t_0|+|\log t_1|+\frac{t_1-t_0}{\eps}\Big)\\
&\le& \frac{2C_{20}}{\eps^2}(1+|\log\eps|+|\log t_0|+|\log t_1|).
\een
Thus, by H\"older's inequality, we have
\ben
\int_{\T^n}|\partial_x\rho_\lb^\eps(x,\hat{t})|dx&\le&\left(\int_{\T^n}\frac{|\partial_x\rho_\lb^\eps(x,\hat{t})|^2}{\rho_\lb^\eps(x,\hat{t})}dx\right)^{1/2}\left(\int_{\T^n}\rho_\lb^\eps(x,\hat{t})dx\right)^{1/2}\\
&\le& \frac{\sqrt{2C_{20}}}{\eps}(1+|\log\eps|+|\log t_0|+|\log t_1|)^{1/2}.
\een
By adjusting the location of $[t_0,t_1]$, we can ensure $\hat t=t$ then conclude the assertion.
\qed

\begin{lem}
\label{lem:est-D2u}
For every $\tau\in(0,1)$ and $t>2\tau$, there exists a constant $C_{22}>0$ such that the following estimate holds
\ben
&&\int_{2\tau}^t \int_{\T^n} (s-\tau)^\alpha e^{-\lb s} |d_x^2 u_\lb^\eps|^2\rho_\lb^\eps(x,s) dxds \\
&\leq& C_{22}\Big(t^\alpha e^{-\lb t}\frac{1+|\log\eps|+|\log t|}{\eps}+\tau^\alpha \frac{1+|\log\eps|+|\log\tau|}{\eps}\\
&&+\lb^{1-\alpha}+\tau^{-1+\alpha}+\lb^{-1-\alpha}\Big)
\een
\end{lem}
\proof
Differentiating \eqref{eq:hj-vis} twice w.r.t. $x$, we obtain
\ben
\lb z_\lb^\eps+\langle \partial_p H(x,d_xu_\lb^\eps), d_xz_\lb^\eps\rangle+{\rm tr}(d_x^2 u_\lb^\eps\partial_{pp}H(x,d_xu_\lb^\eps)d_x^2 u_\lb^\eps)\\
+2{\rm tr}(\partial_{xp}H(x,d_xu_{\lb}^\eps)d_x^2u_\lb^\eps)+\partial_{xx}H(x,d_xu_\lb^\eps)=\eps\Delta z_\lb^\eps,
\een
which further implies 
\be
\label{eq:z}
\lb z_\lb^\eps+\langle \partial_p H(x,d_xu_\lb^\eps), d_xz_\lb^\eps\rangle+\theta|d_x^2 u_\lb^\eps|^2+2{\rm tr}(\partial_{xp}H(x,d_xu_{\lb}^\eps)d_x^2u_\lb^\eps)\\
+\partial_{xx}H(x,d_xu_\lb^\eps)\le \eps\Delta z_\lb^\eps \nonumber
\ee
by the uniformly convexity of $H$, where $z_\lb^\eps:=\Delta u_\lb^\eps$.
Applying \eqref{eq:fp-2} and \eqref{eq:z} to obtain that
\ben
&&\frac{d}{dt}\int_{\T^n} (t-\tau)^\alpha e^{-\lb t}z_\lb^\eps(x) \rho_\lb^\eps(x,t)dx\\
&&=\int_{\T^n} e^{-\lb t}z_\lb^\eps(x)\Big(\alpha(t-\tau)^{\alpha-1} \rho_\lb^\eps(x,t)-\lb(t-\tau)^\alpha \rho_\lb^\eps(x,t)+ (t-\tau)^\alpha \partial_t\rho_\lb^\eps(x,t)\Big) dx\\
&&\ge \int_{\T^n} \alpha(t-\tau)^{\alpha-1}e^{-\lb t}z_\lb^\eps\rho_\lb^\eps+(t-\tau)^\alpha e^{-\lb t} z_\lb^\eps\left({\rm div}(\partial_pH(x,d_xu_\lb^\eps))\rho_\lb^\eps+\eps\Delta\rho_\lb^\eps\right)dx\\
&&+\int_{\T^n}(t-\tau)^\alpha e^{-\lb t} \Big(\langle \partial_p H(x,d_xu_\lb^\eps), d_xz_\lb^\eps\rangle+\theta|d_x^2 u_\lb^\eps|^2-\eps\Delta z_\lb^\eps\\
&&+2{\rm tr}(\partial_{xp}H(x,d_xu_{\lb}^\eps)d_x^2u_\lb^\eps)+\partial_{xx}H(x,d_xu_\lb^\eps)\Big)\rho_\lb^\eps(x,t) dx\\
&&=\int_{\T^n} \alpha(t-\tau)^{\alpha-1}e^{-\lb t}z_\lb^\eps(x)\rho_\lb^\eps(x,t)+\theta(t-\tau)^\alpha e^{-\lb t}|d_x^2 u_\lb^\eps|^2\rho_\lb^\eps(x,t)\\
&&+(t-\tau)^\alpha e^{-\lb t}\Big(2{\rm tr}(\partial_{xp}H(x,d_xu_{\lb}^\eps)d_x^2u_\lb^\eps)+\partial_{xx}H(x,d_xu_\lb^\eps)\Big)\rho_\lb^\eps(x,t) dx,
\een
where the integration by parts formula is used in the last equality. Integrating it from $2\tau$ to $t$, we have
\be
\label{lemA2:est-1}
&&\theta\int_{2\tau}^t \int_{\T^n}(s-\tau)^\alpha e^{-\lb s} |d_x^2 u_\lb^\eps|^2\rho_\lb^\eps(x,s) dxds\\
&&\le \int_{\T^n} (t-\tau)^\alpha e^{-\lb t} z_\lb^\eps(x)\rho_\lb^\eps(x,t) dx-\int_{\T^n} \tau^\alpha e^{-2\lb\tau} z_\lb^\eps(x) \rho_\lb^\eps(x,\tau) dx \nonumber\\
&&-\int_{2\tau}^t\int_{\T^n} \alpha (s-\tau)^{\alpha-1}e^{-\lb s} z_\lb^\eps(x)\rho_\lb^\eps(x,s)dxds \nonumber \\
&&-\int_{2\tau}^t\int_{\T^n} (s-\tau)^\alpha e^{-\lb s}\left(2{\rm tr}(\partial_{xp}H(x,d_xu_{\lb}^\eps)d_x^2u_\lb^\eps)+\partial_{xx}H(x,d_xu_\lb^\eps)\right) \rho_\lb^\eps(x,s) dxds\nonumber
\ee
where $\theta>0$ is a constant related with the convexity of $H$, defined in \eqref{eq:convexity}. Now we estimate the right hand side of the previous inequality. For the first two terms on the right hand side of \eqref{lemA2:est-1}, we can make use of integration by parts. Due to the uniform Lipschitz property of $u_\lb^\eps$ and Lemma \ref{lem:est-Drho}, we get
\ben
\left| \int_{\T^n} (t-\tau)^\alpha e^{-\lb t} z_\lb^\eps(x)\rho_\lb^\eps(x,t)dx\right| &\le& (t-\tau)^\alpha e^{-\lb t}\int_{\T^n} |d_x u_\lb^\eps(x)||\partial_x\rho_\lb^\eps(x,t)|dx \\
&\le& C_1C_{21}(t-\tau)^\alpha e^{-\lb t}\frac{1+|\log\eps|+|\log t|}{\eps}
\een
and
\ben
\left|\int_{\T^n} \tau^\alpha e^{-2\lb\tau} z_\lb^\eps(x) \rho_\lb^\eps(x,\tau)dx\right| &\le& \tau^\alpha \int_{\T^n}|d_xu_\lb^\eps(x)||\partial_x\rho_\lb^\eps(x,\tau)|dx \\
&\le& C_1C_{21}\tau^\alpha \frac{1+|\log\eps|+|\log\tau|}{\eps}.
\een
For the third term on the right hand side of \eqref{lemA2:est-1}, by using Young's inequality, we have
\ben
&&\left|\int_{2\tau}^t\int_{\T^n} \alpha(s-\tau)^{\alpha-1}e^{-\lb s} z_\lb^\eps(x)\rho_\lb^\eps(x,s)dxds\right| \\
&\le& \int_{2\tau}^t\int_{\T^n} \frac{\theta}{3n}(s-\tau)^\alpha e^{-\lb s}|\Delta u_\lb^\eps(x)|^2\rho_\lb^\eps(x,s)+\frac{3n}{4\theta}(s-\tau)^{\alpha-2}e^{-\lb s}\rho_\lb^\eps(x,s)dxds\\
&\le& \frac{\theta}{3}\int_{2\tau}^t (s-\tau)^\alpha e^{-\lb s}|d_x^2u_\lb^\eps(x)|^2\rho_\lb^\eps(x,s)dxds+\frac{3n}{4\theta}\int_{2\tau}^t (s-\tau)^{\alpha-2}e^{-\lb s}ds.
\een
Finally, for the last integral on the right hand side of \eqref{lemA2:est-1}, there exists a positive constant $C_{23}$ such that
\ben
& &\left|\int_{2\tau}^t\int_{\T^n} (s-\tau)^\alpha e^{-\lb s}\Big(2\partial_{xp}H(x,d_xu_{\lb}^\eps)d_x^2u_\lb^\eps+\partial_{xx}H(x,d_xu_\lb^\eps)\Big) \rho_\lb^\eps(x,s) dxds\right|\\
&\le& \int_{2\tau}^t\int_{\T^n} (s-\tau)^\alpha e^{-\lb s} \Big( \frac{3}{\theta}|\partial_{xp}H(x,d_xu_{\lb}^\eps)|^2 +\frac{\theta}{3}|d_x^2 u_\lb^\eps|^2+|\partial_{xx}H(x,d_xu_\lb^\eps)|\Big)\rho_\lb^\eps dxds\\
&\le&  \frac{\theta}{3}\int_{2\tau}^t\int_{\T^n} (s-\tau)^\alpha e^{-\lb s}|d_x^2 u_\lb^\eps|^2 \rho_\lb^\eps(x,s)dxds+ C_{23}\int_{2\tau}^t (s-\tau)^\alpha e^{-\lb s} ds.
\een
Combining the these estimates together and using  Lemma \ref{lemA4} below, we get
\ben
&& \frac{\theta}{3}\int_{2\tau}^t (s-\tau)^\alpha e^{-\lb s}|d_x^2u_\lb^\eps(x)|^2\rho_\lb^\eps(x,s)dxds\\
 &&\le C_1C_{21}(t-\tau)^\alpha e^{-\lb t}\frac{1+|\log\eps|+|\log t|}{\eps}+C_1C_{21}\tau^\alpha \frac{1+|\log\eps|+|\log\tau|}{\eps}\\
 &&+\frac{3n}{4\theta}\Big(\lb^{1-\alpha}+\frac{1}{1-\alpha}\tau^{-1+\alpha}\Big)+C_{23}\lb^{-1-\alpha}\Big(1+\frac{1}{1-\alpha}\Big).
\een
which finishes the proof.
\qed

\begin{lem}
\label{lemA4}
For constants $\tau,\alpha\in(0,1)$ and $\lb t>1$, we have
\beq\label{lemA4:est-1}
\int_{2\tau}^t (s-\tau)^{\alpha} e^{-\lb s}ds \le \lb^{-1-\alpha}\Big(1+\frac{1}{1-\alpha}\Big),
\eeq
\beq\label{lemA4:est-2}
\int_{2\tau}^t (s-\tau)^{-\alpha} e^{-\lb s}ds \le \lb^{-1+\alpha}\Big(1+\frac{1}{1-\alpha}\Big),
\eeq
\beq\label{lemA4:est-3}
\int_{2\tau}^t (s-\tau)^{\alpha-2} e^{-\lb s}ds \le \lb^{1-\alpha}+\frac{1}{1-\alpha}\tau^{-1+\alpha}.
\eeq
\end{lem}
\proof
For $a\in(0,1)$, note that $s^a\le 1$ when $s\le 1$ and $e^{-s}\le s^{-2}$, we have
\ben
\int_{t_0}^{t_1} s^a e^{-s}ds \le \int_{t_0}^1 e^{-s}ds+\int_1^{t_1} s^{a-2}ds=e^{-t_0}-e^{-1}+\frac{1}{a-1}(t_1^{a-1}-1)\le 1+\frac{1}{1-a}.
\een
For $a<0$, note that $e^{-s}\le 1$ and $s^a\le1$ when $s\ge1$, we have
\ben
\int_{t_0}^{t_1} s^a e^{-s}ds \le \int_{t_0}^1 s^a ds+\int_1^{t_1} e^{-s}ds=\frac{1}{a+1}(1-t_0^{a+1})+e^{-1}-e^{-t_1}.
\een
Hence, we obtain the following estimate holds for every $0<t_0<1<t_1$ and $a\in(0,1)$
\be
\int_{t_0}^{t_1} s^a e^{-s}ds\le \left\{
\begin{aligned}
&1+\frac{1}{1-a}, &a\in(0,1),\\
&1+\frac{1}{1+a},  &a\in(-1,0),\\
&1-\frac{1}{1+a}t_0^{a+1}, &a\in(-\infty,-1).
\end{aligned}
\right.
\ee
By using variable substitution, we  get
\ben
\int_{2\tau}^t (s-\tau)^{\alpha} e^{-\lb s}ds=e^{-\lb \tau}\int_\tau^{t-\tau} s^\alpha e^{-\lb s}ds=e^{-\lb\tau}\lb^{-1-\alpha}\int_{\lb\tau}^{\lb(t-\tau)} s^\alpha e^{-s}ds.
\een
\eqref{lemA4:est-1} is therefore given as $\alpha\in(0,1)$. The estimates \eqref{lemA4:est-2} and \eqref{lemA4:est-2} can be obtained similarly by changing $\alpha$ to $-\alpha$ and $\alpha-2$ respectively.
\qed

\begin{thm}\label{thm:3}
For any Tonelli Hamiltonian $H$ and a constant $\alpha\in(1/3,1)$, there exists a positive constant $C_{25}$ uniform w.r.t. $\eps,\lb\in(0,1]$ such that
\[
\|u_\lb^\eps(x)-u_\lb(x)\|\le C_{25}\Big(\frac{\eps}{\lb}\Big)^{\frac{1-\alpha}{2}}+h.o.t.
\]
\end{thm}
\proof
Applying \eqref{eq:hj-partialu} and \eqref{eq:fp-2}, we have
\ben
\frac{d}{dt}\int_{\T^n} e^{-\lb t}v_\lb^\eps(x)\rho_\lb^\eps(x,t)dx&=&\int_{\T^n}e^{-\lb t}[-\lb v_\lb^\eps(x)\rho_\lb^\eps(x,t)+v_\lb^\eps(x)\partial_t\rho_\lb^\eps(x,t)]dx\\
&=&\int_{\T^n} e^{-\lb t}\Big[\left(\langle\partial_p H(x,d_xu_\lb^\eps),d_xv_\lb^\eps\rangle-\Delta u_\lb^\eps-\eps\Delta v_\lb^\eps\right)\rho_\lb^\eps\\
& &+v_\lb^\eps({\rm div}(\partial_p H(x,d_xu_\lb^\eps))\rho_\lb^\eps+\eps\Delta\rho_\lb^\eps)\Big]dx\\
&=&\int_{\T^n}-e^{-\lb t}\Delta u_\lb^\eps(x)\rho_\lb^\eps(x,t) dx,
\een
and integrating it from $0$ to $t$ to obtain that
\be\label{eq:final}
v_\lb^\eps(x)&=&\int_{\T^n} e^{-\lb t}v_\lb^\eps(x)\rho_\lb^\eps(x,t)dx+\int_0^t\int_{\T^n}e^{-\lb s}\Delta u_\lb^\eps(x)\rho_\lb^\eps(x,s)dxds\nonumber\\
&=&\int_{\T^n} e^{-\lb t}v_\lb^\eps(x)\rho_\lb^\eps(x,t)dx+\int_0^{2\tau}\int_{\T^n}e^{-\lb s}\Delta u_\lb^\eps(x)\rho_\lb^\eps(x,s)dxds\nonumber\\
&+&\int_{2\tau}^t\int_{\T^n}e^{-\lb s}\Delta u_\lb^\eps(x)\rho_\lb^\eps(x,s)dxds.
\ee
Due to Proposition \ref{prop:osc-eps}, the first integral of \eqref{eq:final} can be estimated by 
\[
\left|\int_{\T^n} e^{-\lb t}v_\lb^\eps(x)\rho_\lb^\eps(x,t)dx\right|\le \frac{1}{\lb\sqrt{\eps}}e^{-\lb t}.
\]
For the second integral of \eqref{eq:final}, due to \eqref{eq:est-Drho2} and the uniform Lipschitz property of $u_\lb^\eps$ we obtain that
\ben
\left|\int_0^{2\tau}\int_{\T^n}e^{-\lb s}\Delta u_\lb^\eps(x)\rho_\lb^\eps(x,s)dxds\right| &\le& \int_0^{2\tau}\int_{\T^n}e^{-\lb s} |d_xu_\lb^\eps(x)||\partial_x\rho_\lb^\eps(x,s)|dxds\\
&\le& 2C_1C_{21} \tau \frac{1+|\log\eps|+|\log 2\tau|}{\eps}.
\een
Finally, by using Lemma \ref{lem:est-D2u} and the H\"older inequality, we estimate the third integral of \eqref{eq:final}:
\ben
&&\left|\int_{2\tau}^t\int_{\T^n}e^{-\lb s}\Delta u_\lb^\eps(x)\rho_\lb^\eps(x,s)dxds\right| \le \sqrt{n}\int_{2\tau}^t\int_{\T^n} e^{-\lb s}|d_x^2u_\lb^\eps(x)|\rho_\lb^\eps(x,s)dxds\\
&&\le\sqrt{n} \left(\int_{2\tau}^t\int_{\T^n}(s-\tau)^\alpha e^{-\lb s}|d_x^2 u_\lb^\eps(x)|^2\rho_\lb^\eps(x,s)dxds\right)^{1/2}\left(\int_{2\tau}^t(s-\tau)^{-\alpha}e^{-\lb s}ds\right)^{1/2}\\
&&\le \sqrt{C_{22}n}\Big(t^\alpha e^{-\lb t}\frac{1+|\log\eps|+|\log t|}{\eps}+\tau^\alpha \frac{|\log\eps|+|\log \tau|+1}{\eps}\\
&&+\lb^{1-\alpha}+\tau^{-1+\alpha}+\lb^{-1-\alpha}\Big)^{1/2} \Big(1+\frac{1}{1-\alpha}\Big)^{1/2}\lb^{(-1+\alpha)/2}\\
&&\le \sqrt{C_{22}n\Big(1+\frac{1}{1-\alpha}\Big)}\Big((\lb^{-1+\alpha}t^\alpha e^{-\lb t})^{1/2}\frac{(1+|\log t|)^{1/2}+|\log\eps|^{1/2}}{\eps^{1/2}}\\
&&+(\lb^{-1+\alpha}\tau^\alpha )^{1/2}\frac{(1+|\log \tau|)^{1/2}+|\log\eps|^{1/2}}{\eps^{1/2}}+1+(\lb\tau)^{(-1+\alpha)/2}+\lb^{-1}\Big)
\een
where we used the estimate \eqref{lemA4:est-2}.

Combining all these estimates together, we know that there exists a constant $C_{24}>0$ such that
\ben
|v_\lb^\eps(x)|&\le& C_{24}\Big(\frac{e^{-\lb t}}{\lb\sqrt{\eps}}+\frac{\tau}{\eps}(1+|\log\eps|+|\log\tau|)+1+(\lb\tau)^{(-1+\alpha)/2}+\lb^{-1}\\
&+&(\lb^{-1}(\lb t)^\alpha e^{-\lb t})^{1/2}\frac{(1+|\log t|)^{1/2}+|\log\eps|^{1/2}}{\eps^{1/2}}\\
&+&(\lb^{-1+\alpha}\tau^\alpha )^{1/2}\frac{(1+|\log \tau|)^{1/2}+|\log\eps|^{1/2}}{\eps^{1/2}}\Big).
\een
Consequently, by the definition of $v_\lb^\eps$ we have
\ben
&&|u_\lb^\eps(x)-u_\lb(x)|\le \int_0^\eps |v_\lb^{\eps'}(x)|d\eps'\\
&&\le C_{24}\Big( \frac{e^{-\lb t}\sqrt{\eps}}{\lb}+\tau(1+|\log\tau|+|\log\eps|)|\log\eps|+\eps+\frac{\eps}{\lb}+\frac{\eps}{(\lb\tau)^{(1-\alpha)/2}}\\
&&+\Big(\frac{\eps}{\lb}\Big)^{1/2}\Big((\lb t)^\alpha e^{-\lb t}(1+|\log t|)\Big)^{1/2}+(\lb^{-1}(\lb t)^\alpha e^{-\lb t})^{1/2}(\eps+|\log\eps|^2)\\
&&+\Big(\frac{\eps}{\lb}\Big)^{1/2}(\lb \tau)^{\alpha/2}(1+|\log\tau|)^{1/2}+(\lb^{-1+\alpha}\tau^\alpha )^{1/2}(\eps+|\log\eps|^2)\Big).
\een
If $\tau=\eps^\sigma$ and $t=\eps^{-\varrho}/\lb$, then there exists a constant $C_{25}>0$ such that
\ben
& &\|u_\lb^\eps(x)-u_\lb(x)\|\\
&\le& C_{25}\Big( \frac{\eps^{\varrho+{1}/{2}}}{\lb}+\eps^\sigma(1+|\log\eps|)|\log\eps|+\eps+\frac{\eps}{\lb}+\Big(\frac{\eps}{\lb}\Big)^{\frac{1-\alpha}{2}}\eps^{\frac{1+\alpha}{2}-\frac{1-\alpha}{2}\sigma}\\
&+&\Big(\frac{\eps}{\lb}\Big)^{\frac{1}{2}}\eps^{\frac{\varrho(1-\alpha)}{2}}(1+|\log\eps|+|\log\lb|)^{\frac{1}{2}}+\Big(\frac{\eps}{\lb}\Big)^{\frac{1}{2}}\eps^{\frac{\varrho(1-\alpha)-1}{2}}(\eps+|\log\eps|^2)\\
&+&\Big(\frac{\eps}{\lb}\Big)^{\frac{1}{2}}\lb^{\frac{\alpha}{2}}\eps^{\frac{\alpha \sigma}{2}}(1+|\log\eps|)^{1/2}+\Big(\frac{\eps}{\lb}\Big)^{\frac{1-\alpha}{2}}\eps^{\frac{(\sigma+1)\alpha-1}{2}}(\eps+|\log\eps|^2)\Big).
\een 
Hence, as long as 
\beq
\frac{1-\alpha}{\alpha}<\sigma<\frac{1+\alpha}{1-\alpha},\quad \varrho>\frac{1}{1-\alpha}
\eeq
which demands $\alpha\in(1/3,1)$, we can get 
\[
\|u_\lb^\eps(x)-u_\lb(x)\|\le C_{25}\Big(\frac{\eps}{\lb}\Big)^{\frac{1-\alpha}{2}}+h.o.t.
\]
that finished our proof.
\qed
\begin{rmk}
In previous proof, we actually sacrifice the convergence rate of $\|u_\lb^\eps-u_\lb\|$ to remedy the modulus $\eps(\lb)$ to be nearly optimal. That's why the power $\alpha$ couldn't be $1$. We can use the example $\eps(\lb)=\lb^{1+\kappa}$ with $\kappa>0$ to see the loss on the convergence rate: Due to Theorem \ref{thm:low-b}, $\|u_\lb^\eps-u_\lb\|\lesssim \lb^\kappa$, whereas Theorem \ref{thm:3} implies $\|u_\lb^\eps-u_\lb\|\lesssim \lb^{\frac{(1-\alpha)\kappa}{2}}$. As we can see, the latter estimate is much rougher than the former one.
\end{rmk}
\bigskip

\noindent{\it Proof of Theorem \ref{thm:1}:} Due to Theorem \ref{thm:3}, we have 
\[
\|u_\lb^{\eps(\lb)}-u\|\leq \|u_\lb^{\eps(\lb)}-u_\lb\|+\|u_\lb-u\|\leq  C_{25}\Big(\frac{\eps(
\lb)}{\lb}\Big)^{\frac{1-\alpha}{2}}+\|u_\lb-u\|+h.o.t.
\]
for some constant $\alpha\in(0,1/3)$. That implies the convergence of $u_\lb^{\eps(\lb)}$ to $u$ as $\lb\rightarrow 0_+$ once the modulus $\eps(\lb)\leq o(\lb)$.\qed
\\

\section{Examples on the optimal convergence rate of $\|u_\lb^\eps-u_\lb\|$}\label{s5}
In this section, we present an example to show the optimality of the convergence rate of $\|u_\lb^\eps-u_\lb\|$ which we obtained in previous sections. Since obtaining the definite expression of the solution $u_\lb^\eps$ for equation \eqref{eq:hj-vis} is usually impossible, we have to resort to numerical experiments for that.

\begin{ex}[Stochastic pendulum]
Consider the Hamiltonian of one dimensional simple pendulum, i.e., $H(x,p):=\frac{1}{2}|p|^2+\cos x-1$ with $(x,p)\in (\R\slash[0,2\pi])\times\R$. The viscous Hamilton--Jacobi equation can be written as follows
\be
\label{eq:hj-vis-sp}
\lb u_\lb^\eps(x)+\frac{1}{2}|d u_\lb^\eps(x)|^2+\cos x-1=\eps \Delta u_\lb^\eps(x),  \quad x\in \R\slash[0,2\pi].
\ee
Without loss of generality, we take $\eps(\lb)=\lb^{1+\alpha}$ with $\alpha\in(0,1)$. We also need the discount Hamilton--Jacobi equation
\be
\label{eq:hj-dis-sp}
\lambda u_\lb(x)+\frac{1}{2}|du_\lb(x)|^2+\cos x-1=0, \quad x\in\R\slash[0,2\pi].
\ee
Due to the symmetry of \eqref{eq:hj-vis-sp}, $u_\lb(0)\equiv0$ and $d u_\lb(0)=du_\lb^\eps(0)=du_\lb^\eps(\pi)$ for all $\lb,\eps\in(0,1]$. Therefore, 
 \[
  \|u_\lb^\eps(x)-u_\lb(x)\|\geq  \|u_\lb^\eps(0)-u_\lb(0)\|
 \]
 of which the right hand side presents as an effective lower bound. In the following, we will numerically verify that 
  \[
 \|u_\lb^\eps(0)-u_\lb(0)\|\ge C\Big(\frac{\eps}{\lb}+\eps|\log\eps|\Big)
 \]
 for some positive constants $C$.

\begin{figure}[htbp]
  \centering
  \subfigure[$\alpha=0.2$]{
    \includegraphics[width=0.7\textwidth]{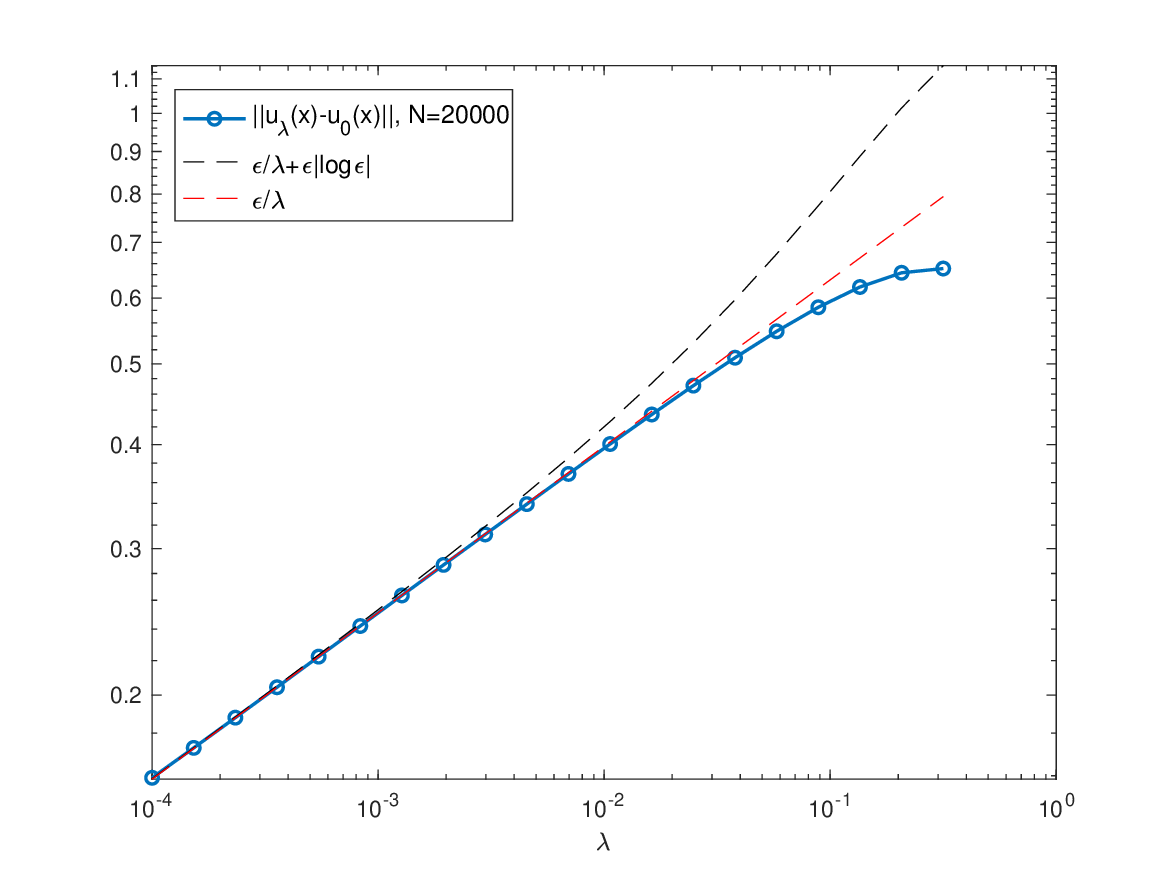}
    \label{fig:subfig1}
  }
  
  \subfigure[$\alpha=0.6$]{
    \includegraphics[width=0.7\textwidth]{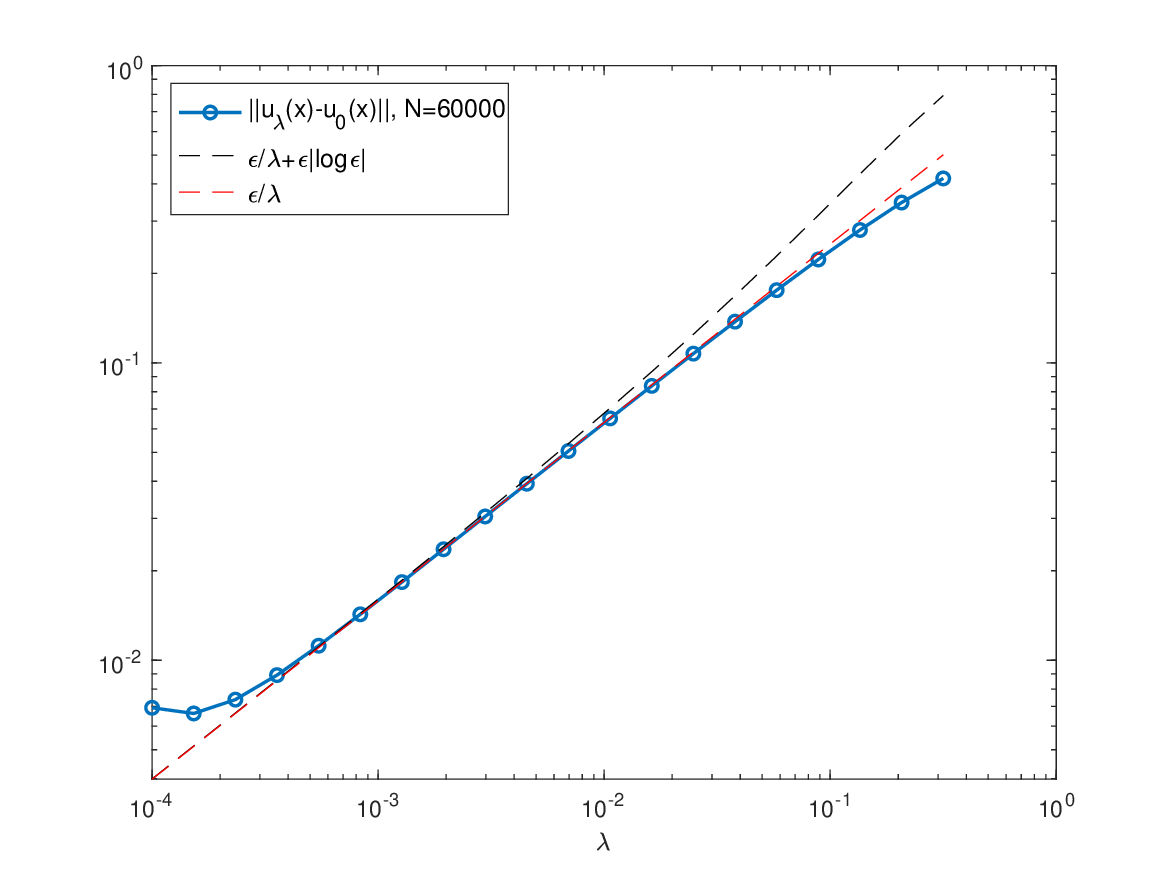}
    \label{fig:subfig2}
  }
  \caption{The convergence rate of  $u_\lb^\eps$ to $u_\lb$ with $\eps=\lb^{1+\alpha}$ as $\lb\to 0_+$. }
  \label{fig:sp}
\end{figure}

Due to the symmetry of equation \eqref{eq:hj-vis-sp}, 
 it suffices to analyze the behavior of $u_\lb^\eps(x)$ and $u_\lb(x)$ on the interval $[0,\pi]$. As we know, 
\[
d_xu_\lb(x)= \sqrt{1-\cos x-\lb u_\lb(x)},  \quad u_\lb(0)=0,
\]
which can be numerically solved by the Runge-Kutta method. However, for \eqref{eq:hj-vis-sp} only a Neumann boundary condition $d u_\lb^\eps(0)=du_\lb^\eps(\pi)=0$ is given, so we employed the numerical methods of finite difference and Newtonian iterations. Notice that 
\[
\lb u_\lb^\eps(0)=\eps\Delta u_\lb^\eps(0),
\]
so it suffices to prove that $|\Dt u_\lb^\eps(0)|\geq C>0$ for all $\lb\in(0,1]$.
Specifically, we partition the interval $[0,\pi]$ into $N$ meshgrid with uniform mesh size $h=\pi/N$. 
Then, the viscous Hamilton--Jacobi equation \eqref{eq:hj-vis-sp} can be discretized as follows
\[
F_j(u)=\lb u_j+\frac{1}{2}\Big(\frac{u_{j+1}-u_{j-1}}{2h}\Big)^2+\cos x_j-1-\eps\frac{u_{j+1}-2u_j+u_{j-1}}{h^2}.
\]
Moreover, the Neumann boundary condition can be embedded in the iterative process by setting additional grid points $u_{-1}$ and $u_{N+1}$. That is 
\[
u'(0)=\frac{u_1-u_{-1}}{2h}=0, \quad u'(\pi)=\frac{u_{N+1}-u_{N-1}}{2h}=0,
\]
with $u_{-1}=u_1$ and $u_{N+1}=u_{N-1}$ (see link \href{https://github.com/zibowangmath/vis\_dis\_HJ}{https://github.com/zibowangmath/vis\_dis\_HJ} for the code).
\end{ex}


\end{document}